\documentstyle{article}

\newcommand{\N}{{\rm I\kern-.5ex N}}
\newcommand{\Z}{{\sf \vrule height 1.55ex depth-1.2ex width.03em\kern-.11em Z
\kern-.9ex Z\kern-.11em\vrule height 0.3ex depth0ex width.03em}}
\newcommand{\Q}{{\rm\kern.2ex\vrule height1.55ex depth-.05ex width.03em\kern-.7ex Q}}
\newcommand{\R}{{\rm I\kern-.5ex R}}
\newcommand{\Rvar}{{\rm I\kern-.5ex R}}
\newcommand{\C}{{\rm\kern.3ex\vrule height1.55ex depth-.05ex width.03em\kern-.7ex C}}
\newcommand{\Cvar}{{\, \rm\kern.3ex\vrule height1.1ex depth-.05ex width.03em\kern-.7ex C}}

\newcommand{\spat}{\hspace{4ex}}

\newcommand{\flip}{ \chi }

\newcommand{\ahh}{\hat{A} \hspace{-.55ex}\hat{\rule{0ex}{2.0ex}}\hspace{.55ex}}

\newcommand{\ah}{\hat{A}}
\newcommand{\deh}{\hat{\Delta}}
\newcommand{\psih}{\hat{\psi}}

\newcommand{\nab}{\nabla}

\newcommand{\od}{\odot}
\newcommand{\ot}{\otimes}

\newcommand{\la}{\Lambda}

\newcommand{\om}{\omega}
\newcommand{\io}{\iota}
\newcommand{\vfi}{\varphi}
\newcommand{\vep}{\varepsilon}

\newcommand{\al}{\alpha}
\newcommand{\be}{\beta}

\newcommand{\sde}{\delta}
\newcommand{\sdeh}{\hat{\sde}}
\newcommand{\de}{\Delta}

\newcommand{\th}{\theta}

\newcommand{\si}{\sigma}

\newcommand{\Mfiu}{{\cal M}_{\vfi_u}}

\newcommand{\text}[1]{\mbox{#1}}

\newcommand{\cst}{\text{C}$\hspace{0.1mm}^*$}
\newcommand{\im}{\text{\rm Im}\,}

\newcommand{\lan}{\langle}
\newcommand{\ran}{\rangle}

\newcommand{\qed}{\ \hfill \rule{2mm}{2mm}}
\newenvironment{demo}{\medskip\noindent\bf Proof :\ \  \rm}{\qed\bigskip\par }

\newtheorem{definition}{Definition}[section]
\newtheorem{proposition}[definition]{Proposition}
\newtheorem{lemma}[definition]{Lemma}
\newtheorem{corollary}[definition]{Corollary}
\newtheorem{remark}[definition]{Remark}
\newtheorem{theorem}[definition]{Theorem}

\newtheorem{result}[definition]{Result}

\begin{document}
\begin{center}
\huge\bf The analytic structure of an algebraic quantum group

\end{center}

\bigskip

\begin{center}
\rm J. Kustermans  \footnote{Research Assistant of the
National Fund for Scientific Research (Belgium)}

Institut for Matematik og Datalogi

Odense Universitet

Campusvej 55

5230 Odense M

Denmark

\bigskip

\bf July 1997 \rm
\end{center}
\bigskip

\subsection*{Abstract}
In \cite{VD1}, Van Daele introduced the notion of an algebraic quantum group. We proved in
\cite{Kus} and \cite{JK} that such algebraic quantum groups give rise to \cst-algebraic
quantum groups according to Masuda, Nakagami \& Woronowicz. In this paper, we will pull
down the analytic structure of these \cst-algebraic quantum groups to the algebraic
quantum group.

\bigskip

\section*{Introduction}

Van Daele introduced the notion of an algebraic quantum group in \cite{VD1}. It is
essentially a Multiplier Hopf-$^*$-algebra with a non-zero left invariant functional on
it. He proved that these algebaic quantum groups form a well-behaved category :
\begin{itemize}
\item The left and right invariant functionals are unique up to a scalar.
\item The left and right invariant functionals are faithful.
\item Each algebraic quantum group gives rise to a dual algebraic quantum group.
\item The dual of the dual is isomorphic to the original algebraic quantum group.
\end{itemize}

This category of algebraic quantum groups contains the compact and discrete quantum
groups. It is also closed under the double construction of Drinfel'd so it contains also
non-compact non-discrete quantum groups. Most of the groups and quantum groups will
however not belong to this category.

It is nevertheless worth while to study these algebraic quantum groups :
\begin{itemize}
\item The category contains interesting examples.
\item The theoretical framework is not technically complicated. In a sense, you have only
to worry about essential quantum group problems. Not over \cst-algebraic complications.
\end{itemize}

\bigskip

At the moment, Masuda, Nakagami \& Woronowicz are working on a possible definition of a
\cst-algebraic quantum group. To get an idea of this definition, we refer to
\cite{MasNak}. You can also get a flavour of it in \cite{Kus} and \cite{JK}. This
definition is technically rather involved and \cst-algebraic quantum groups in this scheme
have a rich analytical sructure.

\medskip

We proved in \cite{Kus} and \cite{JK} that an algebraic quantum group with a positive left
invariant functional gives rise to such \cst-algebraic quantum groups. In this paper, we
will pull down this analytic structure to the algebra level. It shows that such algebraic
quantum groups are truely algebraic versions of \cst-algebraic quantum groups according to
Masuda, Nakagami \& Woronowicz.

\medskip

In the first section, we give an overview of the theory of algebraic quantum groups as can
be found in \cite{VD1}. The second and third section are essentially only intended to
assure that an obvious theory of analytic one-parameter groups works fine on this algebra
level.

The most important results can be found in section 4 where we prove that the analytic
objects of the \cst-algebraic quantum groups are of an algebraic nature.

In the last section, we connect the analytic objects of the dual to the analytic objects
of the original algebraic quantum group.

\bigskip

We end this section with some conventions and notations.

\vspace{1mm}

Every algebra in this paper is an associative algebra over the complex numbers (not
necessarily unital). A homomorphism between algebras is by definition a linear
multiplicative mapping. A $^*$-homomorphism between $^*$-algebras is a homomorphism which
preserves the $^*$-operation.

\vspace{1mm}

If $V$ is a vector space, then $L(V)$ denotes the set of linear mappings from $V$ into
$V$.

If $V$, $W$ are two vector spaces, $V \od W$ denotes the algebraic tensor product. The
flip from $V \od W$ to $W \od V$ will be denoted by $\flip$.

We will also use the symbol $\od$ to denote the algebraic tensor product of two linear
mappings.

\vspace{1mm}

We will always use the minimal tensor product between two \cst-algebras and we will use
the symbol $\ot$ for it. This symbol will also be used to denote the completed tensor
product of two mappings which are sufficiently continuous.

\vspace{1mm}

If $z$ is a complex number, then $S(z)$ will denote the following horizontal strip in the
complex plane : $S(z) = \{ \, y \in \C \mid \im y \in [0,\im z] \, \}$.

\bigskip

\section{Algebraic quantum groups} \label{alg}

In this first section, we will introduce the notion of an algebraic quantum group as can
be found in \cite{VD1}. Moreover, we will give an overview of the properties of this
algebraic quantum group. The proofs of these results can be found in the same paper
\cite{VD1}.  We will first introduce some terminology.

\medskip

We call a $^*$-algebra $A$ non-degenerate if and only if we have for
every $a\in A$ that :
$$(\forall b \in A: a b =0) \Rightarrow a=0 \hspace{1.5 cm} \text{ and  } \hspace{1.5cm}
(\forall b \in A: b a =0) \Rightarrow a=0 .$$

For a non-degenerate $^*$-algebra $A$, you can define the multiplier
algebra $M(A)$. This is a unital $^*$-algebra in which $A$ sits as a
selfadjoint ideal (the definition of this multiplier algebra is the
same as in the case of \cst-algebras).

\medskip

If you have two non-degenerate $^*$-algebras $A,B$ and a
multiplicative linear mapping $\pi$ from $A$ to $M(B)$, we call $\pi$
non-degenerate if and only if the vector spaces $\pi(A) B$ and $B
\pi(A)$ are equal to $B$. Such a non-degenerate multiplicative linear
map has a unique multiplicative linear extension to $M(A)$. This
extension will be denoted by the same symbol as the original mapping.
Of course, we have similar definitions and results for
antimultiplicative mappings. If we work in an algebraic setting, we
will always use this form of non-degeneracy as opposed to the non
degeneracy of $^*$-homomorphisms between \cst-algebras!

\medskip

For a linear functional $\om$ on a non-degenerate $^*$-algebra $A$ and any $a \in M(A)$ we
define the linear functionals $\om a$ and $a \om$ on $A$ such that $(a \om)(x) = \om(x a)$
and $(\om a)(x) = \om(a x)$ for every $x \in  A$.

\medskip

You can find some more information about non-degenerate algebras in the appendix of
\cite{VD6}.

\bigskip

Let $\om$ be a linear functional on a $^*$-algebra $A$, then :
\begin{enumerate}
\item $\om$ is called positive if and only if
$\om(a^* a)$ is positive for every $a \in A$.
\item If $\om$ is positive, then $\om$ is called faithful if and only if
for every $a \in A$, we have that $$\om(a^* a)=0 \Rightarrow a=0.$$
\end{enumerate}

\medskip

Consider a positive linear functional $\om$ on a $^*$-algebra $A$. Let $H$ be a
Hilbert-space and $\la$ a linear mapping from $A$ into $H$ such that
\begin{itemize}
\item $\la(A)$ is dense in $H$.
\item We have for all $a,b \in A$ that $\om(a^* a) = \lan \la(a), \la(b) \ran$.
\end{itemize}
Then we call $(H,\la)$ a GNS-pair for $\om$.

\medskip

Such a GNS-pair always exist and it is unique up to a unitary.

\bigskip

We have now gathered the necessary information to understand the following definition

\begin{definition}
Consider a non-degenerate $^*$-algebra $A$ and a non-degenerate $^*$-homomorphism $\de$
from $A$ into $M(A \od A)$ such that
\begin{enumerate}
\item $(\de \od \io)\de = (\io \od \de)\de$.
\item The linear mappings $T_1$, $T_2$ from $A \od A$ into $M(A \od A)$
such that $$T_1(a \ot b) = \de(a)(b \ot 1)
\hspace{1cm} \text{ and  } \hspace{1cm} T_2(a \ot b) = \de(a)(1 \ot b)$$
for all $a,b \in A$, are bijections from $A \od A$ to $A \od A$.
\end{enumerate}
Then we call $(A,\de)$ a Multiplier Hopf$\,^*$-algebra.
\end{definition}

In \cite{VD6}, A. Van Daele proves the existence of a unique non-zero
$^*$-homomorphism $\vep$ from $A$ to $\C$ such that $$(\vep \od
\io)\de = (\io \od \vep)\de =\io \ .$$ Furthermore, he proves the existence of a
unique anti-automorphism $S$ on $A$ such that
$$m(S \od \io)(\de(a)(1 \ot b)) = \vep(a) b \hspace{1cm} \text{ and } \hspace{1cm}
m(\io \od S)((b \ot 1)\de(a)) = \vep(a) b $$
for every $a,b \in A$  (here, $m$ denotes the multiplication map from $A \od A$ to $A$).
As usual, $\vep$ is called the counit and $S$ the antipode of $(A,\de)$.
Moreover, $S(S(a^*)^*) = a$ for all $a \in A$. Also, $\flip(S \od S)\de = \de S$.

\medskip

Let $\om$ be a linear functional on $A$ an $a$ an element in $A$. We define the element
$(\om \od \io)\de(a)$ in $M(A)$ such that
\begin{itemize}
\item $(\om \od \io)(\de(a)) \, \, b = (\om \od \io)(\de(a)(1 \ot b))$ \item
$b \, \, (\om \od \io)(\de(a)) = (\om \od \io)((1 \ot b)\de(a))$
\end{itemize}
for every $b \in A$.

In a similar way, the multiplier $(\io \od \om)\de(a)$ is defined.

\bigskip

Let $\om$ be a linear functional on $A$. We call $\om$ left invariant
(with respect to $(A,\de)$), if and only if $(\io \od \om)\de(a) =
\om(a)\,1$ for every $a \in A$. Right invariance is defined in a
similar way.

\begin{definition}
Consider a Multiplier Hopf$\,^*$-algebra $(A,\de)$ such that there exists a non-zero
positive linear functional $\vfi$ on $A$ which is left invariant.  Then we call $(A,\de)$
an algebraic quantum group.
\end{definition}

For the rest of this paper, we will fix an algebraic quantum group $(A,\de)$ together with
a non-zero left invariant positive linear functional $\vfi$ on it.

An important feature of such an algebraic quantum group is the
faithfulness and uniqueness of left invariant functionals :
\begin{enumerate}
\item Consider a left invariant linear functional $\om$ on $A$, then there exists a unique
      element $c \in \C$ such that $\om = c \, \vfi$.
\item Consider a non-zero left invariant linear functional $\om$ on $A$, then $\om$ is
      faithful.
\end{enumerate}
In particular, $\vfi$ is faithful.

\medskip

A first application of this uniqueness result concerns the antipode :
Because $\vfi S^2$ is left invariant, there exists a unique complex
number $\mu$ such that $\vfi S^2 = \mu \vfi$ (in \cite{VD1}, our $\mu$
is denoted by $\tau$!). It is not so difficult to prove in an
algebraic way that $|\mu|=1$. The question remains open if there
exists an example of an algebraic quantum group (in our sense) with
$\mu \neq 1$.

\medskip

We define the linear functional $\psi = \vfi S$. It is clear that $\psi$ is a non-zero
right invariant linear functional on $A$. However, in general, $\psi$ will not be
positive. In \cite{Kus}, we use  the \cst-algebra approach to prove the existence of a
non-zero positive right invariant linear functional on $A$.

Of course, we have similar faithfulness and uniqueness results about right invariant
linear functionals.

\medskip

In this paper, we will need frequently the following formula :
\begin{equation}
(\io \od \vfi)(\,(1 \ot a)\de(b)\,)
= S(\,(\io \od \vfi)(\de(a)(1 \ot b))\,)   \label{eq1.3}  \end{equation}
for all $a,b \in A$. A proof of this result can be found in proposition
3.11 of \cite{VD1}. It is in fact nothing else but an algebraic form of the strong left
invariance in the definition of Masuda, Nakagami \& Woronowicz.

\medskip

Another non-trivial property about $\vfi$ is the existence of a unique automorphism $\rho$
on $A$ such that $\vfi(a b) = \vfi(b \rho(a))$ for every $a,b \in A$. We call this the
weak KMS-property of $\vfi$ (In \cite{VD1}, our mapping $\rho$ is denoted by $\sigma$!).

This weak KMS-property is crucial to extend $\vfi$ to a weight on the \cst-algebra level.
We have moreover that $\rho(\rho(a^*)^*) = a$ for every $a \in A$.

As usual, there exists a similar object $\rho'$ for the right
invariant functional $\psi$, i.e. $\rho'$ is an automorphism on $A$
such that $\psi(a b) = \psi(b \rho'(a))$ for every $a,b \in
A$.

Using the antipode, we can connect $\rho$ and $\rho'$ via the formula
$S \rho' = \rho S$. Furthermore, we have that $S^2$ commutes with
$\rho$ and $\rho'$. The interplay between $\rho$,$\rho'$ and $\de$ is
given by the following formulas :
$$\de \rho = (S^2 \od \rho)\de \hspace{1cm} \text{ and } \hspace{1cm} \de \rho' =
(\rho' \od S^{-2})\de.$$

\medskip

It is also possible to introduce the modular function of our algebraic quantum group. This
is an invertible element $\sde$ in $M(A)$ such that $$(\vfi \od \io)(\de(a)(1 \ot b)) =
\vfi(a) \, \sde  \, b $$ for every $a,b \in A$.

Concerning the right invariant functional, we have that $$(\io \od \psi)(\de(a)(b \ot
1)) = \psi(a) \, \sde^{-1} \, b$$ for every $a,b \in A$.

This modular function is, like in the classical group case, a one
dimensional (generally unbounded) corepresentation of our algebraic
quantum group :
$$ \de(\sde) = \sde \od \sde \hspace{2cm} \vep(\sde)= 1 \hspace{2cm}
S(\sde) = \sde^{-1} .$$

As in the classical case, we can relate the left invariant functional
to our right invariant functional via the modular function : we have
for every $a \in A$ that
$$ \vfi(S(a)) =  \vfi(a \sde) = \mu \, \vfi(\sde a) .$$
If we apply this equality two times and use the fact that $S(\sde) =
\sde^{-1}$, we get that $ \vfi(S^2(a)) = \vfi(\sde^{-1} a \sde)$ for every $a \in A$.

\medskip

Not surprisingly, we have also that $\rho(\sde) = \rho'(\sde) =
\mu^{-1} \sde$.

\medskip

Another connection between $\rho$ and $\rho'$ is given by the equality
$\rho'(a) = \sde \rho(a) \sde^{-1}$ for all $a \in A$.

\bigskip

We have also a property which says, loosely speaking, that every
element of $A$ has compact support (see e.g. \cite{JK3} for a proof) :

Consider $a_1,\ldots\!,a_n \in A$. Then there exists an element $c$ in $A$ such that
$c \, a_i = a_i \, c = a_i$ for every $i \in \{1,\ldots\!,n\}$.

\bigskip\bigskip

In a last part, we are going to say something about  duality.

We define the subspace $\ah$ of $A'$ as follows :
$$\ah = \{\, \vfi a \mid a \in A \,\} = \{\, a \vfi \mid a \in A \,\}.$$
Like in the theory of Hopf$\,^*$-algebras, we turn $\ah$ into a
non-degenerate $^*$-algebra :
\begin{enumerate}
\item For every $\om_1,\om_2 \in \ah$ and $a \in A$, we have that $(\om_1 \om_2)(a) =
(\om_1 \od \om_2)(\de(a))$.
\item For every $\om \in \ah$ and $a \in A$, we have that $\om^*(a) =
\overline{\om(S(a)^*)}$.
\end{enumerate}
We should remark that a little bit of care has to be taken by defining the product and
the $^*$-operation in this way.

\bigskip

In \cite{JK3}, we proved that $M(\ah)$ can be implemented as a subspace of $A'$ (theorem
2.12 and proposition 7.4)

\begin{proposition}
As a vector space, $M(\ah)$ is equal to
$$\{\, \th \in A' \mid \text{ We have for every } a \in A
\text{ that } (\th \od \io)\de(a) \text{ and } (\io \od \th)\de(a) \text{ belong to }
A \,\} \ . $$
The multiplication and $^*$ operation on $M(\ah)$ are defined in the following way.
\begin{enumerate}
\item We have for every $\th_1,\th_2 \in M(\ah)$ and $a \in A$ that
$(\th_1 \th_2)(a) = \th_1( (\io \od \th_2)\de(a) )$  $= \th_2( (\th_1 \od \io)\de(a) )$.
\item We have for every $\th \in M(\ah)$ and $a \in A$ that $\th^*(a) =
\overline{\th(S(a)^*)}$.
\end{enumerate}
Furthermore, $\vep$ is the unit of $M(\ah)$.
\end{proposition}

\medskip

This implies also that $M(\ah \od \ah)$ can be viewed as a subspace of $(A \od A)'$. Then
we can use this implementation of $M(\ah \od \ah)$ to define a comultiplication $\deh$ on
$\ah$ such that $\deh(\om)(x \ot y) = \om(x y)$ for every $\om \in \ah$ and $x,y \in A$.

\vspace{1mm}

A definition of the comultiplication $\deh$ without the use of such an
embedding can be found in definition $4.4$ of \cite{VD1}.

\medskip

In this way, $(\ah,\deh)$ becomes a Multiplier Hopf$\,^*$-algebra. The counit $\hat{\vep}$
and the antipode $\hat{S}$ are such that
\begin{enumerate}
\item We have for every $\om \in A'$ and $a \in A$ that $\hat{\vep}(\om a)=
\hat{\vep}(a \om )= \om(a)$.
\item We have for every $\om \in M(\ah)$ and $a \in A$ that $\hat{S}(\om)(a) = \om(S(a))$.
\end{enumerate}

For any $a \in A$, we define $\hat{a} = a \vfi \in \ah$. The mapping $A \rightarrow \ah :
a \mapsto \hat{a}$ is a bijection, which is in fact nothing else but the Fourier
transform.

\medskip

Next, we define the linear functional $\psih$ on $\ah$ such that
$\psih(\hat{a}) = \vep(a)$ for every $a \in A$. It is possible to prove that $\psih$
is right invariant.

We have also that $\psih(\hat{a}^* \hat{a}) = \vfi(a^* a)$ for every $a \in A$. This
implies that $\psih$ is a non-zero positive right invariant linear functional on $\ah$.

\medskip

We will also define a linear functional $\hat{\vfi}$ on $\ah$ such that $\hat{\vfi}(\psi
a) = \vep(a)$ for $a \in A$. Then $\hat{\vfi}$ is a non-zero left invariant functional on
$\ah$. It is not clear wether $\hat{\vfi}$ is positive.

\bigskip

From theorem 9.9 of \cite{Kus}, we know  that $(A,\de)$ possesses a non-zero positive
right invariant linear functional.  In a similar way, this functional will give rise to a
non-zero positive left invariant linear functional on $\ah$. This will imply that
$(\ah,\deh)$ is again an algebraic quantum group.

\bigskip

\section{Analytic one-parameter groups}

In this section, we introduce the notion of analytic one-parameter groups on an algebraic
quantum group. It is an obvious variant of the well-known concept of one-parameter groups
on a \cst-algebra. As a consequence, the proofs in this section are essentially the same
as in the \cst-algebra case.

\medskip

These analytic one-parameter groups will be essential in a later section in order to
describe the polar decomposition of the antipode and the modular properties of the Haar
functional on the $^*$-algebra level.

\bigskip

We will fix an algebraic quantum group $(A,\de)$ with a positive left invariant functional
$\vfi$ on it. Let $(H,\la)$ be a GNS-pair of $\vfi$.

\medskip

Let us start of with the definition of an analytic one-parameter group on $(A,\de)$.

\begin{definition}
Consider a function $\al$ from $\C$ into the set of homomorphisms from $A$ into $A$ such
that the following properties hold :
\begin{enumerate}
\item We have for every $t \in \R$ that $\al_t$ is a  $^*$-automorphism on $A$.
\item We have for every $s,t \in \R$ that $\al_{s+t} = \al_s \al_t$.
\item We have for every $t \in \R$ that $\al_t$ is relatively invariant under $\vfi$
\item Consider $a \in A$ and $\om \in \ah$. Then the function
$\C \rightarrow \C : z \mapsto \om(\al_z(a))$ is analytic.
\end{enumerate}
Then we call $\al$ an analytic one-parameter group on $(A,\de)$.
\end{definition}

\bigskip

Except for the last proposition in this section, we will fix an analytic one-parameter
group $\al$ on $(A,\de)$. We will prove the basic properties of $\al$.

\medskip

\begin{result}
Consider $z \in \C$.
Then $\al_z(a)^* = \al_{\overline{z}}(a^*)$ for $a \in A$.
\end{result}
\begin{demo}
Choose $a \in A$.

Take $b \in A$. Then we have two analytic functions
$$\C \rightarrow \C : u \mapsto \overline{\vfi(\al_{\overline{u}}(a^*) \,b^*)}
\hspace{1.5cm} \text{ and } \hspace{1.5cm}
\C \rightarrow \C : u \mapsto \vfi(b \, \al_u(a))$$
These functions are equal on the real axis : we have for all $t \in \R$ that
$$\overline{\vfi(\al_{\overline{t}}(a^*)\,b^*)} = \overline{\vfi(\al_t(a^*) \, b^*)}
= \vfi(b \, \al_t(a))$$
So they must be equal on the whole complex plane.

We have in particular that $\overline{\vfi(\al_{\overline{z}}(a^*)\,b^*)} =
\vfi(b\,\al_z(a))$ which implies that $\vfi(b\, \al_{\overline{z}}(a^*)^*) = \vfi(b \,
\al_z(a))$.

So the faithfulness of $\vfi$ implies that $\al_{\overline{z}}(a^*)^* =  \al_z(a)$.
\end{demo}

\medskip

In the next result, we extend the group character of $\al$ to the whole complex plane.

\begin{result}
Consider $y,z \in \C$. Then $\al_{y+z} = \al_y \, \al_z$.
\end{result}
\begin{demo} Choose $a \in A$.

Take $b \in A$. Fix $s \in \R$ for the moment. Then there exists by assumption a strictly
positive number $M$ such that $\vfi \, \al_s = M \, \vfi$. We have two analytic functions
:
$$\C \rightarrow \C : u \mapsto M \, \vfi(\al_{-s}(b) \al_u(a)) \hspace{1.5cm}
\text{ and } \hspace{1.5cm}
\C \rightarrow \C : u \mapsto  \vfi(b \, \al_{s+u}(a))$$
These functions are equal on the real axis : we have for all $t \in \R$ that
$$M \, \vfi(\al_{-s}(b) \al_t(a)) = \vfi\bigl(\al_s(\al_{-s}(b) \al_t(a))\bigr)
= \vfi\bigl(b \, \al_s(\al_t(a))\bigr) = \vfi(b \, \al_{s+t}(a))$$
So these two functions must be equal on the whole complex plane.

In particular, $M \, \vfi(\al_{-s}(b) \al_z(a)) = \vfi(b \, \al_{s+z}(a))$. Therefore, we
get that
$$\vfi(b \, \al_{s+z}(a)) = M \, \vfi(\al_{-s}(b) \al_z(a)) = \vfi\bigl(\al_s(\al_{-s}(b)
\al_z(a))\bigr) = \vfi\bigl(b \, \al_s(\al_z(a))\bigr) \hspace{1.5cm} \text{(*)}$$

\medskip

Now we have again two analytic functions :
$$\C \rightarrow \C : u \mapsto  \vfi(b \, \al_u(\al_z(a))\bigr) \hspace{1.5cm}
\text{ and } \hspace{1.5cm}
\C \rightarrow \C : u \mapsto  \vfi(b \, \al_{u+z}(a))$$
which by (*) are equal on the real axis. So they must be equal on the whole complex plane.

Hence we get in particular that $\vfi(b\,\al_{y+z}(a)) =
\vfi\bigl(b\,\al_y(\al_z(a))\bigr)$.

\vspace{1mm}

Therefore the faithfulness of $\vfi$ implies that $\al_{y+z}(a) = \al_y(\al_z(a))$.
\end{demo}

\begin{corollary}
Consider $z \in \C$. Then $\al_z$ is an automorphism on $A$ and $(\al_z)^{-1} = \al_{-z}$.
\end{corollary}

\bigskip

We want to use $\al$ to define a positive injective operator in  the GNS-space $H$ of
$\vfi$. In order to do so, we will need the following results.

\begin{result}
There exists a unique strictly positive number $\lambda$ such that
$\vfi \, \al_t = \lambda^t \, \vfi$ for $t \in \R$.
\end{result}
\begin{demo}
By assumption, there exists for every $t \in \R$ a strictly positive element $\lambda_t$
such that $\vfi \, \al_t = \lambda_t \, \vfi$.

It is then clear that $\lambda_0 = 1$ and that $\lambda_s \, \lambda_t = \lambda_{s+t}$
for every $s,t \in \R$.

\vspace{1mm}

Now there exist $a,b \in A$ such that $\vfi(b a) \neq 0$.

We have that the function $\R \rightarrow \C : t \rightarrow
\vfi(\al_{-t}(b)\,a)$ is continuous and non-zero in 0. Hence there exist a
strictly positive number $t_0$ such that $\vfi(\al_{-t}(b) \, a) \neq 0$ for $t \in
[-t_0,t_0]$.

Now we have for every $t \in [-t_0,t_0]$ that $\vfi(b \, \al_t(a)) =
\vfi(\al_t(\al_{-t}(b)\,a)) = \lambda_t \,\vfi(\al_{-t}(b)\,a)$ which implies that
$\lambda_t = \frac{\vfi(b \,\al_t(a))}{\vfi(\al_{-t}(b)\,a)}$. So we see that the function
$[-t_0,t_0] \rightarrow \R^+_0 : t \mapsto \lambda_t$ is continuous.

Therefore the function $\R \rightarrow \R^+_0 : t \mapsto \lambda_t$ is continuous.

We get from this all the existence of a strictly positive number $\lambda$ such that
$\lambda^t = \lambda_t$ for $t \in \R$.  So $\vfi\,\al_t = \lambda^t \, \vfi$ for $t \in
\R$.
\end{demo}

\medskip

The above quality can now be generalized to the whole complex plane.

\begin{result}
We have that $\vfi \, \al_z = \lambda^z \, \vfi$ for every $z \in \C$.
\end{result}
\begin{demo}
Choose $a,b \in A$.
Then we have two analytic functions
$$\C \rightarrow \C : u \mapsto \lambda^u \, \vfi(a \, \al_{-u}(b))
 \hspace{1.5cm} \text{and} \hspace{1.5cm}
\C \rightarrow \C : u \mapsto  \vfi(\al_u(a) \, b)$$
These functions are equal on the real line : we have for all $t \in \R$ that
$$\lambda^t \, \vfi(a \, \al_{-t}(b)) = \vfi(\al_t(a \, \al_{-t}(b)))
= \vfi(\al_t(a) \, b)$$
So they must be equal on the whole complex plane.

We have in particular that $\vfi(\al_z(a)\, b) = \lambda^z \, \vfi(a \, \al_{-z}(b))$,
hence $\vfi\bigl(\al_z(a \, \al_{-z}(b))\bigr) = \lambda^z \, \vfi(a \, \al_{-z}(b))$.

Because $A = A \, \al_{-z}(A)$, we infer from this that $\vfi \, \al_z = \lambda^z \,
\vfi$.
\end{demo}

\medskip

Now we can define a natural positive injective operator in $H$.

\begin{definition}
There exists a unique injective positive operator $P$ in $H$ such that \newline
$P^{it} \la(a) = \lambda^{-\frac{t}{2}} \, \la(\al_t(a))$ for $a \in A$ and $t \in \R$.
\end{definition}
\begin{demo}
We can define a unitary group representation $u$ from $\R$ on $H$ such that $u_t \, \la(a)
= \lambda^{-\frac{t}{2}} \, \la(\al_t(a))$ for every $t \in \R$ and $a \in A$.

Choose $a,b \in A$. Then we have for every $t \in \R$ that
$$\lan u_t \, \la(a) , \la(b) \ran
= \lambda^{-\frac{t}{2}} \, \lan \la(\al_t(a)) , \la(b)
= \lambda^{-\frac{t}{2}} \, \vfi(b^* \al_t(a))$$
This implies that the function $\R \rightarrow \C : \lan u_t \, \la(a) , \la(b) \ran$ is
continuous. Because $u$ is bounded and $\la(A)$ is dense in $H$, we conclude that the
function $\R \rightarrow \C : \lan u_t \, v , w \ran$ is continuous for every $v,w \in H$.

So we see that $u$ is weakly and hence strongly continuous.

By the Stone theorem, there exists a unique injective positive operator $P$ in $H$ such
that $P^{it} = u_t$ for every $t \in \R$.
\end{demo}

\medskip

We want to show that every $P^{iz}$ is defined on $\la(A)$ and that the formula in the
above definition has its obvious generalization to $P^{iz}$. First we will need a lemma
for this.

\begin{lemma}
Consider $a \in A$. Then  the function $\C \rightarrow \C : z \mapsto \lambda^{- \text{\rm
\scriptsize Re}\, z} \, \vfi(\al_z(a)^* \al_z(a))$ is bounded on horizontal strips.
\end{lemma}
\begin{demo}
Take $r \in \R$. We will prove that the function above is bounded on the horizontal strip
$S(r i)$.

We have for every $t \in \R$ that $$\vfi(\al_{ti}(a)^* \al_{ti}(a)) = \vfi(\al_{-ti}(a^*)
\al_{t i}(a)) = \lambda^{-it} \, \vfi(a^* \al_{2 t i}(a))$$ which implies that the
function $\R \rightarrow \R^+ : t \mapsto \vfi(\al_{ti}(a)^* \al_{ti}(a))$ is continuous.

So there exists $M \in \R^+$ such that $\vfi(\al_{ti}(a)^* \al_{ti}(a)) \leq M$ for $t \in
[0,r]$.

Hence we get for every $z \in \C$ that
$$\lambda^{-\text{\rm \scriptsize Re}\, z} \, \vfi(\al_{z}(a)^* \al_{z}(a))
= \lambda^{-\text{\rm \scriptsize Re}\, z} \, \vfi\bigl(\al_{\text{\rm
\scriptsize Re}\, z}(\al_{\text{\rm \scriptsize Im}\, z}(a))^*
\al_{\text{\rm \scriptsize Re}\, z}(\al_{\text{\rm \scriptsize Im}\, z}(a))\bigr)
= \vfi(\al_{\text{\rm \scriptsize Im}\, z}(a)^* \al_{\text{\rm \scriptsize Im}\, z}(a))
\leq M$$
\end{demo}

\medskip

\begin{result}  \label{art2.res1}
Consider  $z \in \C$. Then $\la(A) \subseteq D(P^{i z})$ and $P^{i z} \la(a) =
\lambda^{-\frac{z}{2}} \, \la(\al_z(a))$ for $a \in A$.
\end{result}
\begin{demo} Take $v \in H$ and define the function $f$ from $\C$ into $\C$ such that
$f(u) = \lan \lambda^{-\frac{u}{2}} \, \la(\al_u(a)) , v \ran$ for $u \in \C$.

Define also for every $b \in A$ the function $f_b : \C \rightarrow \C$ such that
$f_b(u) = \lan \lambda^{-\frac{u}{2}} \, \la(\al_u(a)) , \la(b) \ran$ for $u \in \C$.

Then $f_b(u) = \lambda^{-\frac{u}{2}}  \, \vfi(b^* \al_u(a))$ for $u \in \C$ which implies
that $f_b$ is analytic.

\medskip

Choose $y \in \C$ and consider the open unit bal $B$ in $\C$ around $y$.

By the previous lemma, we get the existence of a positive number $M$ such that
$\lambda^{-\text{\rm \scriptsize Re}\, u} \, \vfi(\al_u(a)^* \al_u(a)) \leq M^2$ for $u
\in B$. This implies easily that $\| \lambda^{-\frac{u}{2}} \la(\al_u(a)) \| \leq M$ for
$u \in B$.

Now there exists a sequence $(b_n)_{n=1}^\infty$ in $A$ such that
$(\la(b_n))_{n=1}^\infty$ converges to $v$. We have for every $u \in B$ and $n \in \N$
that
$$|f(u) - f_{b_n}(u)| = |\lan \lambda^{-\frac{u}{2}} \, \la(\al_u(a)) , \la(b_n) \ran
- \lan \lambda^{-\frac{u}{2}} \, \la(\al_u(a)) , v \ran | \leq M \, \|\la(b_n) - v\|$$

This implies that $(f_{b_n})_{n=1}^\infty$ converges uniformly to $f$ on $B$, so $f$ is
analytic on $B$. This implies that $f$ is analytic in $y$.

\vspace{1mm}

So we see that $f$ is analytic on $\C$. It is also clear that
$f(t) = \lan P^{it} \la(a) , v \ran$ for $t \in R$.

\medskip

This implies that $\la(a) \in D(P^{iz})$ and that $P^{i z} \la(a) = \lambda^{-\frac{z}{2}}
\la(\al_z(a))$.
\end{demo}

\begin{result}
Consider $z \in \C$. Then $\la(A)$ is a core for $P^{iz}$.
\end{result}

This follows from the fact that $\la(A)$ is dense in $H$ and that $P^{it} \la(A) \subseteq
\la(A)$ for $t \in \R$ (see e.g. corollary 1.22 of \cite{JK1}).

\bigskip\medskip

In the last proposition of this section, we prove that analytic one-parameter groups are
determined by their value in $i$.

\medskip

\begin{proposition}  \label{art2.prop1}
Consider analytic one-parameter groups $\al$,$\be$ on $(A,\de)$.
If $\al_i = \be_i$, then $\al = \be$.
\end{proposition}
\begin{demo} It is clear that also $\al_{-i} = \be_{-i}$.

Now there exist strictly positive numbers $\lambda$ and $\mu$ such that
$\vfi \, \al_t = \lambda^t \, \vfi$ and $\vfi \, \be_t = \mu^t \, \vfi$ for $t \in \R$.

Then there exist also positive injective operators $P$ and $Q$ in $H$ such that
$P^{it} \la(a) = \lambda^{-\frac{t}{2}} \, \la(\al_t(a))$ and
$Q^{it} \la(a) = \mu^{-\frac{t}{2}} \, \la(\be_t(a))$ for $t \in \R$ and $a \in A$.

We know that $\la(A)$ is a core for both $P$ and $Q$. We have moreover for every $a \in A$
that $$\mu^\frac{i}{2} \, P \la(a) = \mu^\frac{i}{2} \lambda^\frac{i}{2} \la(\al_{-i}(a))
= \lambda^\frac{i}{2} \mu^\frac{i}{2} \la(\be_{-i}(a))
= \lambda^\frac{i}{2} \, Q \la(a)$$

So we get that $\mu^\frac{i}{2} \, P = \lambda^\frac{i}{2} \, Q$. It is clear that
$\frac{\mu^\frac{i}{2}}{\lambda^\frac{i}{2}}$ is a positive number. Because is also has
modulus 1, it must be equal to 1. So $P = Q$.

Hence we get for every $a \in A$ and $z \in \C$ that $\lambda^{-\frac{z}{2}} \la(\al_z(a))
= \mu^{-\frac{z}{2}} \la(\be_z(a))$ which by the faitfulness of $\vfi$ implies that
$\lambda^{-\frac{z}{2}} \, \al_z(a) = \mu^{-\frac{z}{2}} \be_z(a)$.

Now we see that $\frac{\mu}{\lambda} \, \al_2 = \be_2$, so
we get that $\frac{\mu}{\lambda} \, \al_2$ is multiplicative.
Because also $\al_2$ is multiplicative, this implies that $\frac{\mu}{\lambda} = 1$.
Hence $\al = \be$.
\end{demo}

\bigskip

\section{Analytic unitary representations}

This section contains the same basic ideas as the previous one. We only apply them in a
different situation. We will use the results in this section to introduce complex powers
of the modular function of an algebraic quantum group in the next section.

\bigskip

We will again fix an algebraic quantum group $(A,\de)$ with a positive left Haar
functional $\vfi$ on it. Let $(H,\de)$ be a GNS-pair for $\vfi$.

\medskip

Let us start of with a definition.

\begin{definition}
Consider a function $u$ from $\C$ into $M(A)$ such that the following properties hold :
\begin{enumerate}
\item We have for every $t \in \R$ that $u_t$ is  unitary in $M(A)$.
\item We have for every $s,t \in \R$ that $u_{s+t} = u_s u_t$.
\item Consider $\om \in \ah$. Then the function
$\C \rightarrow \C : z \mapsto \om(u_z(a))$ is analytic.
\end{enumerate}
Then we call $u$ an analytic unitary representation on $A$.
\end{definition}

\bigskip

Except for the last result of this section, we will fix an analytic one-parameter
representation $u$ on $A$. We will prove the basic properties of $u$.

\medskip

\begin{result}
We have that $u_z^* = u_{-\overline{z}}$ for every $z \in \C$.
\end{result}
\begin{demo}
Take $a \in A$. Then we have two analytic functions
$$\C \rightarrow \C : y \mapsto \overline{\vfi(u_{-\overline{y}}\,a^*) }
\hspace{1.5cm} \text{ and } \hspace{1.5cm}
\C \rightarrow \C : y \mapsto \vfi(a \, u_y)$$
These two functions are equal on the real line : we have for all $t \in \R$ that
$$\overline{\vfi( u_{-\overline{t}}\,a^*)}  = \overline{\vfi(u_{-t}\, a^*)}
= \vfi(a \, u_{-t}^*) = \vfi(a \,u_t)$$
for all $t \in \R$. So they must be equal on the whole complex plane. We have in
particular that $\overline{\vfi( u_{-\overline{z}} \, a^*)} = \vfi(a \, u_z)$ which
implies that $\vfi(a \, u_z) = \vfi(a \, u_{-\overline{z}}^*)$.

So the faithfulness of $\vfi$ implies that $u_{-\overline{z}}^* =  u_z$.
\end{demo}

\medskip

Now we extend the group character of $u$ to the whole complex plane.

\begin{result}
Consider $y,z \in \C$. Then $u_{y+z} = u_y \, u_z$.
\end{result}
\begin{demo} Choose $a \in A$.

Fix $s \in \R$ for the moment.  We have two analytic functions :
$$\C \rightarrow \C : c \mapsto \vfi((a u_s)\,u_c) \hspace{1.5cm}
\text{ and } \hspace{1.5cm}
\C \rightarrow \C : c \mapsto  \vfi(a \, u_{s+c})$$
which are equal on the real axis. So they  must be equal on the whole complex plane.
In particular, $\vfi(a \, u_s \, u_z) = \vfi(a \, u_{s+z})$. \ \ \ (*)

\medskip

Now we have again two analytic functions :
$$\C \rightarrow \C : c \mapsto  \vfi(a \, u_c \, u_z) \hspace{1.5cm}
\text{ and } \hspace{1.5cm}
\C \rightarrow \C : c \mapsto  \vfi(a \, u_{c+z})$$
which by (*) are equal on the real axis. So they must be equal on the whole complex plane.

Hence we get in particular that $\vfi(a \, u_{y+z}) = \vfi(a \, u_y \, u_z)$.

Therefore the faithfulness of $\vfi$ implies that $u_{y+z} = u_y \, u_z$.
\end{demo}

\begin{corollary}
Consider $z \in \C$. Then $u_z$ is invertible in $M(A)$ and $(u_z)^{-1} = u_{-z}$.
\end{corollary}

\bigskip

We next use $u$ to define a positive injective operator in  the GNS-space $H$ of $\vfi$.

\begin{definition}
There exists a unique injective positive operator $P$ in $H$ such that
$P^{it} \la(a) = \la(u_t \, a)$ for $a \in A$ and $t \in \R$.
\end{definition}
\begin{demo}
We can define a unitary group representation $v$ from $\R$ on $H$ such that
$v_t \, \la(a) =  \la(u_t \, a)$ for every $t \in \R$ and $a \in A$.

Choose $a,b \in A$. Then we have for every $t \in \R$ that
$$\lan v_t \, \la(a) , \la(b) \ran
= \lan \la(u_t\, a) , \la(b) \ran =  \vfi(b^* u_t \, a)$$
This implies that the function $\R \rightarrow \C : \lan v_t \, \la(a) , \la(b) \ran$ is
continuous. Because $v$ is bounded and $\la(A)$ is dense in $H$, we conclude that the
function $\R \rightarrow \C : \lan u_t \, v , w \ran$ is continuous for every $v,w \in H$.

So we see that $u$ is weakly and hence strongly continuous.

By the Stone theorem, there exists a unique injective positive operator $P$ in $H$ such
that $P^{it} = u_t$ for every $t \in \R$.
\end{demo}

\medskip

We want to show that every $P^{iz}$ is defined on $\la(A)$ and that the formula in the
above definition has its obvious generalization to $P^{iz}$. First we will need a lemma
for this.

\begin{lemma}
Consider $a \in A$. Then function $\C \rightarrow \C : z \mapsto \vfi((u_z \, a)^*
(u_z\,a))$ is bounded on horizontal strips.
\end{lemma}
\begin{demo}
Take $r \in \R$. We will prove that the function above is bounded on the horizontal strip
$S(r i)$.

We have for every $t \in \R$ that $$\vfi((u_{ti}\,a)^* (u_{ti}\,a)) = \vfi(a ^* u_{ti}\,
u_{t i} \, a) =  \vfi(a^* u_{2 t i} \, a)$$ which implies that the function
$\R \rightarrow \R^+ : t \mapsto \vfi((u_{ti} \, a)^* (u_{ti} \, a))$ is continuous.

So there exists $M \in \R^+$ such that $\vfi((u_{ti}\, a)^* (u_{ti}\,a)) \leq M$ for $t
\in [0,r]$.

Hence we get for every $z \in \C$ that
\begin{eqnarray*}
\vfi((u_z \, a)^* (u_z \, a)) & = & \vfi\bigl((u_{\text{\rm \scriptsize Re}\, z} \,
u_{i\,\text{\rm \scriptsize Im}\, z} \, a)^* \, u_{\text{\rm \scriptsize Re}\, z} \,
u_{i\,\text{\rm \scriptsize Im}\, z}\, a\bigr)
= \vfi((u_{i\,\text{\rm \scriptsize Im}\, z} \, a)^* \, u_{\text{\rm \scriptsize Re}\,
z}^*\, u_{\text{\rm \scriptsize Re}\, z} \,  (u_{i\,\text{\rm \scriptsize Im}\, z}\,
a))  \\
& = & \vfi((u_{i\,\text{\rm \scriptsize Im}\, z}\, a)^*
(u_{i\,\text{\rm \scriptsize Im}\, z}\,a)) \leq M
\end{eqnarray*}
\end{demo}

\medskip

\begin{result}
Consider  $z \in \C$. Then $\la(A) \subseteq D(P^{i z})$ and $P^{i z} \la(a) = \la(u_z \,
a)$ for $a \in A$.
\end{result}
\begin{demo} Take $v \in H$ and define the function $f$ from $\C$ into $\C$ such that
$f(c) = \lan \la(u_c \,  a) , v \ran$ for $c \in \C$.

Define also for every $b \in A$ the function $f_b : \C \rightarrow \C$ such that
$f_b(c) = \lan \la(u_c \, a) , \la(b) \ran$ for $c \in \C$.

Then $f_b(c) =  \vfi(b^* u_c  \, a)$ for $c \in \C$ which implies that
$f_b$ is analytic.

\medskip

Choose $y \in \C$ and consider the open unit bal $B$ in $\C$ around $y$.

By the previous lemma, we get the existence of a positive number $M$ such that
$\vfi((u_c \, a)^* (u_c \, a)) \leq M^2$ for $c \in B$.
This implies easily that $\| \la(u_c\, a) \|  \leq M$ for $c \in B$

Now there exists a sequence $(b_n)_{n=1}^\infty$ in $A$ such that
$(\la(b_n))_{n=1}^\infty$ converges to $v$. We have for every $c \in B$ and $n \in \N$
that
$$|f(c) - f_{b_n}(c)| = |\lan \la(u_c \, a) , \la(b_n) \ran
- \lan \la(u_c \, a) , v \ran | \leq M \, \|\la(b_n) - v\|$$

This implies that $(f_{b_n})_{n=1}^\infty$ converges uniformly to $f$ on $B$, so $f$ is
analytic on $B$. This implies that $f$ is analytic in $y$.

\vspace{1mm}

So we see that $f$ is analytic on $\C$. It is also clear that
$f(t) = \lan P^{it} \la(a) , v \ran$ for $t \in R$.

\medskip

This implies that $\la(a) \in D(P^{iz})$ and that $P^{i z} \la(a) =
\la(u_z \, a)$.
\end{demo}

\begin{result}
Consider $z \in \C$. Then $\la(A)$ is a core for $P^{iz}$.
\end{result}

This follows from the fact that $\la(A)$ is dense in $H$ and that $P^{it} \la(A) \subseteq
\la(A)$ for $t \in \R$ (see e.g. corollary 1.22 of \cite{JK1}).

\bigskip\medskip

In the last proposition of this section, we prove that analytic unitary representations
are determined by their value in $i$.

\medskip

\begin{result}
Consider analytic unitary representations $u$,$v$ on $A$.
If $u_i = v_i$, then $u = v$.
\end{result}
\begin{demo} It is clear that also $u_{-i} = v_{-i}$.

Now there exist also positive injective operators $P$ and $Q$ in $H$ such that
$P^{it} \la(a) = \la(u_t\,a)$ and
$Q^{it} \la(a) = \la(v_t\,a)$ for $t \in \R$ and $a \in A$.

We know that $\la(A)$ is a core for both $P$ and $Q$. We have moreover for every $a \in A$
that $$P \la(a) = \la(u_{-i}\,a)
= \la(v_{-i}\,(a)) =  Q \la(a)$$

So we get that $P = Q$. Hence we get for every $a \in A$ and $z \in \C$ that
$\la(u_z \, a) = \la(v_z \, a)$ which by the faitfulness of $\vfi$ implies that
$u_z \, a = v_z \, a$.
\end{demo}

\bigskip

\section{The analytic structure of an algebraic quantum group}

In \cite{Kus} and \cite{JK}, we proved that every algebraic quantum group gives rise to a
reduced and universal \cst-algebraic quantum group in the sense of Masuda, Nakagami and
Woronowicz. These \cst-algebraic quantum groups have a very rich (complicated?) analytic
structure. We show in this section that this analytic structure can be completely pulled
down to the algebraic level.

\bigskip

Consider an algebraic quantum group $(A,\de)$ and let $\vfi$ be a positive left Haar
functional of $(A,\de)$. Let $(H,\la)$ be a GNS-pair for $\vfi$.

\medskip

In \cite{JK}, we constructed a universal \cst-algebraic quantum group $(A_u,\de_u)$ out of
$(A,\de)$.

We denote the canonical embedding of $A$ into $A_u$ by $\pi_u$. So $\pi_u$ is an injective
$^*$-homomorphism into $A_u$ such that $\pi_u(A)$ is dense in $A_u$. We have also for
every $a \in A$ and $x \in A \od A$ that
$$\de_u(\pi_u(a)) \, (\pi_u \od \pi_u)(x) = (\pi_u \od \pi_u)(\de(a) \, x)
\hspace{0.5cm} \text{ and } \hspace{0.5cm}
(\pi_u \od \pi_u)(x) \, \de_u(\pi_u(a))  = (\pi_u \od \pi_u)(x \, \de(a))$$

\medskip

As a rule, we will give objects associated with $A_u$ a subscript \lq u\rq. So the left
Haar weight on $(A_u,\de_u)$ will be denoted by $\vfi_u$ and its modular group by $\si_u$.

\medskip

The scaling group of $(A_u,\de_u)$ will be denoted by $\tau_u$, the anti-unitary antipode
by $R_u$.

The modular group of the right Haar weight $\vfi_u \, R_u$ will be denoted by $\si_u'$.

\medskip

We will also need the co-unit on the universal \cst-algebraic quantum group $(A_u,\de_u)$.
Recall that we have a co-unit $\vep$ on the algebraic quantum group $(A,\de)$. This gives
rise to a non-zero $^*$-homomorphism $\vep_u$ from $A_u$ into $\C$ such that $\vep_u \,
\pi_u = \vep$. Then we have that $(\vep_u \ot \io)\de_u = (\io \ot \vep_u)\de_u = \io$.

The usefulness of the antipode is the only reason to prefer the universal \cst-algebraic
quantum group over the reduced one in this section.

\bigskip

In \cite{Kus}, we constructed a reduced \cst-algebraic quantum group $(A_r,\de_r)$ out of
$(A,\de)$.

We denote the canonical embedding of $A$ into $A_r$ by $\pi_r$. So $\pi_r$ is an injective
$^*$-homomorphism into $A_r$ such that $\pi_r(A)$ is dense in $A_r$.

Objects associated to $(A_r,\de_r)$ will get the subscript \lq r\rq. So will we denote the
left Haar weight on $(A_r,\de_r)$  by $\vfi_r$ and its modular group $\si_r$.

\medskip

By the universality of $A_u$, there exists a unique $^*$-homomorphism $\pi$ from $A_u$
into $A_r$ such that $\pi \, \pi_u = \pi_r$.

This mapping $\pi$ connects al objects of $A_u$ and $A_r$. So have we for instance that
$\vfi_r \, \pi = \vfi_u$ and $\pi \, (\si_u)_t = (\si_r)_t \, \pi$ for $t \in \R$.

\bigskip\medskip

The first object that we will pull down to the algebra level is the modular function
$\sde_u$ of the \cst-algebraic quantum group. So $\sde_u$ is a strictly positive element
affiliated with $A_u$. The main result to pull down $\sde_u$ to the algebra level is
already proven in proposition 8.12 of \cite{Kus} and proposition 12.11 of \cite{JK}.

\medskip

By proposition 12.11 of \cite{JK}, we get for every $a \in A$ and $z \in \C$ that
$\pi_u(a)$ belongs to $D(\sde_u^{iz})$ and that $\sde_u^{iz} \, \pi_u(a)$ belongs again to
$\pi_u(A)$. We will use this fact to define any complex power of $\sde$ on the algebra
level. The definition that we introduce in this paper is much better than the ad hoc
definitions in \cite{Kus} and \cite{JK} because we work here within the framework of the
algebraic quantum group $(A,\de)$.

\medskip

\begin{proposition} \label{art1.prop7}
There exist a unique analytic unitary representation $w$ on $A$ such that
$w_{-i} = \sde$. We define $\sde^z = w_{-iz}$ for $z \in \C$.
\end{proposition}
\begin{demo}
We know already that $\sde_u^{i y} \, \pi_u(A) \subseteq \pi_u(A)$ for $y \in \C$.

It is then not very difficult to prove the existence of a mapping $w$ from $\C$ into
$M(A)$ such that $\pi_u(w_y \, a) = \sde_u^{i y} \, \pi_u(a)$ for $y \in \C$ and $a \in A$
and such that $w_y^* = w_{-\overline{y}}$ for $y \in \C$.

Then it easy to check that
\begin{enumerate}
\item $w_0 = 1$
\item We have for every $s,t \in \R$ that $w_{s+t} = w_s \, w_t$.
\end{enumerate}

\medskip

Take $\om \in A$. Then there exists $a,b,c \in A$ such that $\om = ab \vfi c^*$.

We know that the $\pi_u(a)$ is analytic with respect to $\sde_u$. This means that the
mapping $\C \rightarrow A_u : y \mapsto \sde_u^{iy} \, \pi_u(a)$ is analytic. This gives
us that the mapping $\C \rightarrow \C : y \mapsto \vfi_u(\pi_u(c)^* \, (\sde_u^{i y} \,
\pi_u(a)) \, \pi_u(b))$ is analytic.

But we have for every $y \in \C$ that
$$\vfi_u(\pi_u(c)^* \, (\sde_u^{i y} \, \pi_u(a)) \, \pi_u(b))
= \vfi_u(\pi_u(c)^* \, \pi_u(w_y \, a) \, \pi_u(b))
= \vfi_u(\pi_u(c^* \, w_y \, a b)) = \om(w_y)$$ so the function
$\C \mapsto \C : y \mapsto \om(w_y)$ is analytic.

\medskip

So we get by definition that $w$ is an analytic one-parameter representation on $A$.
Proposition 12.11 of \cite{JK} gives us that $w_{-i} = \sde$.
\end{demo}

\medskip

Notice that we have also proven the following characterization for powers of $\sde$.

\begin{proposition}
We have for every $z \in \C$ and $a \in A$ that
$\pi_u(\sde^z \, a) = \sde_u^z \, \pi_u(a)$ for $a \in A$.
\end{proposition}

So we have in particular that $\pi_u(\sde^{it} \, a) = \sde_u^{it} \, \pi_u(a)$ and
$\pi_u(a \, \sde^{it}) = \pi_u(a) \, \sde_u^{it}$ for $a \in A$ and $t \in \R$.

\medskip

We will quickly recall the basic calculation rules for powers of $\sde$ (which are
immediate consequences of the results of the previous section).

\begin{proposition}
We have the following calculation rules :
\begin{enumerate}
\item $\sde^0 = 1$ and $\sde^1 = \sde$
\item We have for $y,z \in \C$ that $\sde^{y+z} = \sde^y \, \sde^z$
\item We have for $z \in \C$ that $(\sde^{z})^* = \sde^{\overline{z}}$
\item Consider $z \in \C$. Then $\sde^z$ is invertible in $M(A)$ and
$(\sde^{z})^{-1} = \sde^{-z}$.
\item Consider $t \in \R$. Then $\sde^{it}$ is a unitary element in $M(A)$.
\item Consider $t \in \R$. Then $\sde^{t}$ is a selfadjoint element in $M(A)$.
\end{enumerate}
\end{proposition}

\begin{remark} \rm Let $t$ be a real number. Then $\sde^t$ is  positive in the
sense that $\sde^t = \sde^\frac{t}{2} \, \sde^\frac{t}{2}$ where $\sde^\frac{t}{2}$ is
self adjoint.
\end{remark}

\medskip

We have also the following analyticity property.

\begin{result}
We have for every $\om \in \ah$ that the function $\C \rightarrow \C : z \rightarrow
\om(\sde^z)$ is analytic.
\end{result}

\medskip\medskip

By proposition 12.2 of \cite{JK}, we know that $\pi(\sde_u) = \sde_r$. Hence we get that
$\pi(\sde_u^z) = \sde_r^z$ for $z \in \C$. This implies then easily the next proposition.

\begin{proposition}
We have for every $z \in \C$ and $a \in A$ that
$\pi_r(\sde^z \, a) = \sde_r^z \, \pi_r(a)$ for $a \in A$.
\end{proposition}

So we have in particular that $\pi_r(\sde^{it} \, a) = \sde_r^{it} \, \pi_u(a)$ and
$\pi_r(a \, \sde^{it}) = \pi_r(a) \, \sde_u^{it}$ for $a \in A$ and $t \in \R$.

\bigskip\bigskip

In the next part of this section, we are going to pull down the other objects $\si_u$,
$\si_u'$, $\tau_u$ and $R_u$ to the algebra level.

\medskip

By result 12.8 of \cite{JK}, we get for every $t \in \R$ that
$$(\si_u')_t(\pi_u(A)) = (\si_u)_t(\sde_u^{it} \, \pi_u(A) \, \sde_u^{-it})
= (\si_u)_t(\pi_u(\sde^{it} \, A \, \sde^{-it})) = (\si_u)_t(\pi_u(A))$$

\medskip

\begin{proposition} \label{art1.prop1}
Consider $t \in \R$. Then we have the equalities
$$(\si_u)_t(\pi_u(A)) =  \pi_u(A)  \hspace{2cm} (\si_u')_t(\pi_u(A)) =  \pi_u(A)
\hspace{2cm} (\tau_u)_t(\pi_u(A)) = \pi_u(A) $$
\end{proposition}
\begin{demo} Take $s \in \R$.

Choose $a \in A$. Because $\vep_u \neq 0$, we get that $\vep_u (\si_u)_s \neq 0$. Hence
the density of $\pi_u(A)$ in $A_u$ implies the existence of $b \in A$ such that
$\vep_u\bigl((\si_u)_s(\pi_u(b))\bigr) = 1$.

Now there exist $p_1,\ldots\!,p_n , q_1,\ldots\!,q_n \in A$ such that
$\de(a)(1 \ot b) = \sum_{i=1}^n p_i \ot q_i$.

Then $$\de_u(\pi_u(a))\, (1 \ot \pi_u(b)) = \sum_{i=1}^n \pi_u(p_i) \ot
\pi_u(q_i)\hspace{1.5cm}\text{(*)}$$

By result 9.4 of \cite{JK}, we know that $((\tau_u)_s \ot (\si_u)_s)\de_u = \de_u \,
(\si_u)_s$. So if we apply $(\tau_u)_s \ot (\si_u)_s$ to the equation above, we get that
$$\de_u\bigl((\si_u)_s(\pi_u(a))\bigr)\,\bigl(1 \ot (\si_u)_s(\pi_u(b))\bigr) =
\sum_{i=1}^n (\tau_u)_s(\pi_u(p_i)) \ot (\si_u)_s(\pi_u(q_i))$$
If we now apply $\io \ot \vep_u$ to the equation above and use the fact that $(\io \ot
\vep_u)\de_u = \io$, we see that
$$(\si_u)_s(\pi_u(a)) = \sum_{i=1}^n  (\tau_u)_s(\pi_u(p_i))
\ \vep_u\bigl((\si_u)_s(\pi_u(q_i))\bigr)$$
which implies that
$$((\tau_u)_{-s} \, (\si_u)_s)(\pi_u(a)) = \sum_{i=1}^n  \pi_u(p_i)
\ \vep_u\bigl((\si_u)_s(\pi_u(q_i))\bigr)$$

So we see that $((\tau_u)_{-s}\,(\si_u)_s)(\pi_u(A)) \subseteq \pi_u(A)$. \ \ \ (a)

\medskip

We have also the equality $((\si_u')_s \ot (\tau_u)_{-s})\de_u = \de_u \,(\si_u')_s$
(result 10.16 of \cite{JK}). This gives us in a similar way that $((\tau_u)_s \,
(\si_u')_s)(\pi_u(A)) \subseteq \pi_u(A)$. So by the remarks before this proposition, we
get that $((\tau_u)_s \, (\si_u)_s)(\pi_u(A)) \subseteq \pi_u(A)$ which implies by that
$((\si_u)_s \, (\tau_u)_s)(\pi_u(A)) \subseteq \pi_u(A)$. \ \ \ (b)

\medskip

Therefore (a) and (b) imply that
$(\si_u)_{2 s} (\pi_u(A))  \subseteq \pi_u(A)$.

\medskip\medskip

As a consequence, we have for every $r \in R$ that
$(\si_u)_r(\pi_u(A)) \subseteq \pi_u(A)$. This implies for every $r \in \R$ that
$(\si_u)_r(\pi_u(A)) = \pi_u(A)$.

Combining this with equality (a), we get that $(\tau_u)_r(\pi_u(A)) \subseteq \pi_u(A)$
for $r \in \R$.
\end{demo}

\bigskip

We want now to prove a generalization of this result for complex parameters. The idea
behind this proof is completely the same but we have to be a little bit careful in this
case.

\medskip

Recall from \cite{JK} that every element of $\pi_u(A)$ is analytic with respect to
$\si_u$, $\si_u'$ and $\tau_u$.

\bigskip

First we prove a little lemma.

\begin{lemma} \label{art1.lem1}
Consider $z \in \C$. Then $(\si_u)_z(\pi_u(A))$ is dense in $A_u$.
\end{lemma}
\begin{demo}
Choose $\om \in A_u^*$ such that $\om = 0$ on $(\si_u)_z(\pi_u(A))$.

Take $x \in \pi_u(A)$.
Define the analytic function $f$ from $\C$ into $\,\C$ such that
$f(c) = \om((\si_u)_c(x))$ for $c \in \C$.

We have by the previous proposition for every $t \in \R$ that $(\si_u)_t(x)$ belongs to
$\pi_u(A)$, which implies that $\om\bigl((\si_u)_z((\si_u)_t(x))\bigr) = 0$. Hence $f(z+t)
= \om((\si_u)_{z+t}(x)) = \om\bigl((\si_u)_z((\si_u)_t(x))\bigr) = 0$ for $t \in \R$.

This implies that $f = 0$. In particular, $\om(x) = f(0) = 0$.

So the density of $\pi_u(A)$ in $A_u$ implies that $\om = 0$.

\medskip

Therefore Hahn-Banach implies that $(\si_u)_z(\pi_u(A))$ is dense in $A_u$.
\end{demo}

\bigskip

\begin{lemma} \label{art1.lem2}
Consider $z \in \C$ and $a \in A$. Then $(\si_u')_z(\pi_u(a)) = (\si_u)_z(\pi_u(\sde^{i z}
\,a \, \sde^{-i z}))$.
\end{lemma}
\begin{demo}
We introduce the one-parameter representations $L$, $R$, $\th$ on $A_u$ such that
$$L_t(x) = \sde_u^{it} \, x \hspace{2cm} R_t(x) = x \, \sde_u^{-i t} \hspace{2cm}
\th_t(x) = \sde_u^{it} \, x \, \sde_u^{-it}$$
for $x \in A_u$ and $t \in \R$.

Then $L$ and $R$ commute and $\th_t = L_t \, R_t$ for $t \in \R$. This implies that $L_z
\, R_z \subseteq \th_z$ (see e.g. proposition 3.9 of \cite{JK1}).

We know that $\pi(a) \in D(R_z)$ and that $R_z(\pi_u(a)) = \pi_u(a \, \sde_u^{-iz})$. This
implies that $R_z(\pi_u(a)) \in D(L_z)$ and that $L_z(R_z(\pi_u(a))) = \pi_u(
\sde_u^{iz}\,a \, \sde_u^{-iz})$. So we get that $\pi_u(a) \in D(\th_z)$ and that
$\th_z(\pi_u(a)) = \pi_u( \sde_u^{iz}\,a \, \sde_u^{-iz})$.

\medskip

We have moreover that $\si_u$ and $\th$ commute and that $(\si_u')_t = (\si_u)_t \, \th_t$
for all $t \in \R$. This implies again that $(\si_u)_z \, \th_z \subseteq (\si_u')_z$.
Therefore $(\si_u')_z(\pi_u(a)) = (\si_u)_z\bigl(\th_z(\pi_u(a))\bigr) =
(\si_u)_z(\pi_u(\sde^{iz}\, a \, \sde^{-iz}))$.
\end{demo}

\medskip

This lemma implies immediately the following one.

\begin{lemma}
Consider $z \in \C$. Then $(\si_u)_z(\pi_u(A)) = (\si_u')_z(\pi_u(A))$.
\end{lemma}

\bigskip

Now we are in a position to prove a generalization of proposition \ref{art1.prop1}.

\begin{proposition} \label{art1.prop2}
Consider $z \in \C$. Then we have the equalities
$$(\si_u)_z(\pi_u(A)) =  \pi_u(A)  \hspace{2cm} (\si_u')_z(\pi_u(A)) =  \pi_u(A)
\hspace{2cm} (\tau_u)_z(\pi_u(A)) = \pi_u(A) $$
\end{proposition}
\begin{demo}
Take $y \in \C$.

Consider elements $a_1,\ldots\!,a_m , b_1,\ldots\!,b_m \in A$ and $p_1,\ldots\!,p_n ,
q_1,\ldots\!,q_n \in A$ such that  the equality \newline $\sum_{i=1}^m \de(a_i)(1 \ot b_i)
= \sum_{j=1}^n p_j \ot q_j$ holds. Then $$\sum_{i=1}^m \de_u(\pi_u(a_i))\,(1 \ot
\pi_u(b_i)) = \sum_{j=1}^n \pi_u(p_j) \ot \pi_u(q_j)$$

We have now the following two analytic functions :
$$ \C \rightarrow A_u \ot A_u : c \mapsto \sum_{i=1}^m \de_u\bigl((\si_u)_c(\pi_u(a_i))
\bigr) \,\bigl(1 \ot (\si_u)_c(\pi_u(b_i))\bigr) $$
and
$$ \C \rightarrow A_u \ot A_u : c \mapsto
\sum_{j=1}^n (\tau_u)_c(\pi_u(p_j)) \ot (\si_u)_c(\pi_u(q_j))$$
These two functions are equal on the real line :  we have for every $t \in \R$ that
\begin{eqnarray*}
& & \sum_{i=1}^m \de_u\bigl((\si_u)_t(\pi_u(a_i))\bigr) \, \bigl(1 \ot
(\si_u)_t(\pi_u(b_i))\bigr)
= \sum_{i=1}^m ((\tau_u)_t \ot (\si_u)_t)\bigl(\de_u(\pi_u(a_i))\bigr)
\, \bigl(1 \ot (\si_u)_t(\pi_u(b_i))\bigr) \\
& & \spat = \sum_{i=1}^m ((\tau_u)_t \ot (\si_u)_t)\bigl(\de_u(\pi_u(a_i))
(1 \ot \pi_u(b_i))\bigr)
= \sum_{j=1}^n (\tau_u)_t(\pi_u(p_j)) \ot (\si_u)_t(\pi_u(q_j))
\end{eqnarray*}
So these functions are equal on the whole complex plane. Hence,
$$\sum_{i=1}^m \de_u\bigl((\si_u)_y(\pi_u(a_i))\bigr)\,\bigl(1 \ot
(\si_u)_y(\pi_u(b_i))\bigr)
= \sum_{j=1}^n (\tau_u)_y(\pi_u(p_j)) \ot (\si_u)_y(\pi_u(q_j))$$

Applying $\io \ot \vep_u$ to this equation and using the fact that $(\io \ot \vep_u)\de_u
= \io$ gives us now that
$$\sum_{i=1}^m (\si_u)_y(\pi_u(a_i)) \ \vep_u\bigr((\si_u)_y(\pi_u(b_i))\bigr)
= \sum_{j=1}^n (\tau_u)_y(\pi_u(p_j)) \ \vep_u\bigl((\si_u)_y(\pi_u(q_j))\bigr)
\hspace{1.5cm} \text{(*)}$$

Because $(\si_u)_y(\pi_u(A))$ is dense in $A_u$ (lemma \ref{art1.lem1}) and $\vep_u \neq
0$, we get the existence of $e \in A$ such that $\vep_u\bigl((\si_u)_y(\pi_u(e))\bigr) =
1$. Combining this with the fact that $\de(A)(1 \ot A) = A \od A$ and equation (*), it is
not very difficult to see that $(\si_u)_y(\pi_u(A)) = (\tau_u)_y(\pi_u(A))$. \ \ \ (a)

\medskip

Using the equation $((\si_u')_t \ot (\tau_u)_{-t})\de_u = \de_u \, (\si_u')_t$ for every
$t \in \R$, we get in a similar way that  $(\si_u')_{-y}(\pi_u(A)) = (\tau_u)_y(\pi_u(A))$
which implies that $(\si_u)_{-y}(\pi_u(A)) = (\tau_u)_y(\pi_u(A))$ by the previous lemma.

\medskip

From (a) and this last equality, we get that $(\si_u)_y(\pi_u(A)) =
(\si_u)_{-y}(\pi_u(A))$.

\vspace{1mm}

Take $v \in \pi_u(A)$. By the previous equality, there exists  an element $w \in \pi_u(A)$
such that $(\si_u)_y(v) = (\si_u)_{-y}(w)$. Because $(\si_u)_{-y} = ((\si_u)_{y})^{-1}$,
this implies that $(\si_u)_y(v) \in D((\si_u)_y)$ and $(\si_u)_y((\si_u)_y(v)) = w$. Hence
$(\si_u)_{2 y}(v) = w \in \pi_u(A)$.

\vspace{1mm}

Consequently, $(\si_u)_{2 y}(\pi_u(A)) \subseteq \pi_u(A)$.

\medskip

So we get that $(\si_u)_z(\pi_u(A)) \subseteq \pi_u(A)$.
Combining this with (a), we get moreover that
$(\tau_u)_z(\pi_u(A)) \subseteq \pi_u(A)$.
\end{demo}

\medskip

\begin{corollary} \label{art1.cor1}
We have that $R_u(\pi_u(A)) = \pi_u(A)$.
\end{corollary}
\begin{demo}
We know by theorem 9.18 of \cite{JK} that $R_u((\tau_u)_{-\frac{i}{2}}(\pi_u(a))) =
\pi_u(S(a))$ for every $a \in A$. Therefore the previous proposition implies that
$R_u(\pi_u(A)) = \pi_u(A)$.
\end{demo}

\bigskip\medskip

These results imply that al the objects associated to the \cst-algebraic quantum group
$(A_u,\de_u)$ can be pulled down to the algebraic level. The first (and typical) result is
contained in the next proposition.

\begin{proposition} \label{art1.prop3}
There exists a unique analytic one-parameter group $\si$ on $A$ such that
$\si_{-i} = \rho$.

We have moreover that $\pi_u(\si_z(a)) = (\si_u)_z(\pi_u(a))$ for every $a \in A$ and $z
\in \C$.
\end{proposition}
\begin{demo}
By proposition \ref{art1.prop2}, we know that $(\si_u)_z(\pi_u(A)) = \pi_u(A)$ for every
$z \in \C$.

This implies that we can define a mapping $\si$ from $\C$ into the set of mappings from
$A$ into $A$ such that $\pi_u(\si_z(a)) = (\si_u)_z(\pi_u(a))$ for every $a \in A$ and $z
\in \C$. Then we get immediately that
\begin{enumerate}
\item We have for every $z \in \C$ that $\al_z$ is a homomorphism on $A$.
\item We have for every $t \in \R$ that $\al_t$ is a $^*$-automorphism on $A$.
\item We have for every $s,t \in \R$ that $\al_{s+t} = \al_s \, \al_t$.
\end{enumerate}

\medskip

Choose $z \in A$. Then $\pi_u(a) \in D((\si_u)_z) \cap \Mfiu$ and
$(\si_u)_z(\pi_u(a)) = \pi_u(\si_z(a)) \in \Mfiu$.

Because $\vfi_u$ is invariant under $\si_u$, this implies by proposition 2.14 of
\cite{JK2} that $$\vfi_u(\pi_u(a)) = \vfi_u\bigl((\si_u)_z(\pi_u(a))\bigr) =
\vfi_u\bigl(\pi_u(\si_z(a))\bigr)$$ which implies that $\vfi(a) = \vfi(\si_z(a))$.

\medskip

Take $a \in A$ and $\om \in \ah$. Then there exist $b,c \in A$ such that $\om = c \vfi
b^*$.

We know by proposition 8.8 of \cite{JK} that $\pi_u(a)$ is analytic with respect to
$\si_u$. This implies that the function $\C \rightarrow A_u : z \mapsto
(\si_u)_z(\pi_u(a))$ is analytic. Hence the function $\C \rightarrow \C : z \mapsto
\vfi_u( \pi_u(b)^* \, (\si_u)_z(\pi_u(a)) \, \pi_u(c))$ is analytic.

But we have for every $z \in \C$ that
\begin{eqnarray*}
& & \vfi_u( \pi_u(b)^*\,(\si_u)_z(\pi_u(a))\,\pi_u(c)) = \vfi_u(\pi_u(b)^*
\pi_u(\si_z(a)) \pi_u(c)) \\
& & \spat = \vfi_u(\pi_u(b^* \si_z(a) c)) = \vfi(b^* \si_z(a) c) = \om(\si_z(a))
\end{eqnarray*}

So the function $\C \rightarrow \C : z \mapsto \om(\si_z(a))$ is analytic.

\medskip

Hence we get by definition that $\si$ is an analytic one-parameter group on $A$.

Using proposition 8.8 of \cite{JK}, we see that $\si_{-i} = \rho$.
\end{demo}

We  call $\si$ the modular group of the left Haar weight.

\medskip

We have also proven the following result.

\begin{result}
We have for every $z \in \C$ that $\vfi \, \si_z = \vfi$.
\end{result}

\bigskip

There exists a canonical GNS-construction $(H,\la_u,\pi)$ for the weight $\vfi_u$ such
that $\pi_u(A)$ is a core for $\la_u$ and $\la_u(\pi_u(a)) = \la(a)$ for $a \in A$ (see
definition 10.2 of \cite{JK} and theorem 10.6 of \cite{JK}).

Denote the modular operator of $\vfi_u$ by $\nab$ and the modular conjugation of $\vfi_u$
by $J$ (both with respect to this GNS-construction).

Then these objects are completely characterized by the following two results.

\begin{result}
Consider $z \in \C$. Then $\la(A)$ is a core for $\nab^{iz}$ and
that $\nab^{iz} \la(a) = \la(\si_z(a))$ for ever $a \in A$.
\end{result}
\begin{demo} We have for every $t \in \R$ and $a \in A$ that
$$\nab^{it} \la(a) = \nab^{it} \la_u(\pi_u(a)) = \la_u((\si_u)_t(\pi_u(a))) =
\la_u(\pi_u(\si_t(a))) = \la(\si_t(a))$$

The result follows now from result \ref{art2.res1}.
\end{demo}

\begin{corollary}
We have for every $a \in A$ that $J \la(a) = \la(\si_\frac{i}{2}(a)^*)$
\end{corollary}

\medskip\medskip

By the remarks after definition 8.4 of \cite{JK}, we know that $\pi \, (\si_u)_t =
(\si_r)_t \, \pi$ for $t \in \R$. Consequently, $\pi \, (\si_u)_z \subseteq (\si_r)_z \,
\pi$ for $z \in \C$. So we get easily the next result.

\begin{proposition}
We have for every $z \in \C$ and $a \in A$ that $(\si_r)_z(\pi_r(a)) = \pi_r(\si_z(a))$.
\end{proposition}

\bigskip\bigskip

Completely similar to proposition \ref{art1.prop3}, we have the next results. Remember
that we have proven in corollary 7.3 of \cite{Kus} the existence of a unique strictly
positive number $\nu$ such that $\vfi_r \, (\tau_r)_t = \nu^t \, \vfi_r$ for $t \in \R$.
Proposition 9.19 of \cite{Kus}, proposition 8.17 of \cite{Kus} and proposition 8.20 of
\cite{JK2} imply that $\vfi_r \, (\si_r')_t = \nu^{-t} \, \vfi_r$ for $t \in \R$.

\vspace{1mm}

This gives us also that $\vfi_r \, (\tau_u)_t = \nu^t \, \vfi_u$ and that $\vfi_u \,
(\si_u)_z = \nu^{-t} \, \vfi_u$ for every $t \in \R$.

\bigskip

\begin{proposition}
There exists a unique analytic one-parameter group $\si'$ on $A$ such that
$\si'_{-i} = \rho'$.

We have moreover that $\pi_u(\si_z'(a)) = (\si_u')_z(\pi_u(a))$ for every $a \in A$ and $z
\in \C$.
\end{proposition}

We call $\si'$ the modular function of the right Haar functional.

\medskip

\begin{result}
We have for every $z \in \C$ that $\vfi \, \si_z' = \nu^{-z} \, \vfi$.
\end{result}

\begin{proposition}
We have for every $a \in A$ and $z \in \C$ that $\pi_r(\si_z'(a)) = (\si_r')_z(\pi_r(a))$.
\end{proposition}

\bigskip

In the next part, we pull down the polar decomposition of the antipode to the algebra
level. The proof of the next proposition is completely similar to the proof of proposition
\ref{art1.prop3}.

\begin{proposition}
There exists a unique analytic one-parameter group $\tau$ on $A$ such that
$\tau_{-i} = S^2$. We have moreover that
$\pi_u(\tau_z(a)) = (\tau_u)_z(\pi_u(a))$ for every $a \in A$ and $z \in \C$.
\end{proposition}

We call $\tau$ the scaling group of $(A,\de)$.

\medskip

\begin{result} \label{art1.res3}
We have for every $z \in \C$ that $\vfi \, \tau_z = \nu^z \, \vfi$.
\end{result}

\begin{proposition}
We have for every $a \in A$ and $z \in \C$ that $\pi_r(\tau_z(a)) = (\tau_r)_z(\pi_r(a))$.
\end{proposition}

\bigskip

Now it is the turn of the anti-unitary antipode $R_u$ to be pulled down. This is possible
thanks to corollary \ref{art1.cor1}

\begin{definition}
We define the mapping $R$ from $A$ into $A$ such that
$\pi_u(R(a)) = R_u(\pi_u(a))$ for $a \in A$.
Then $R$ is a $^*$-anti-automorphism on $A$ such that $R^2 = \io$.
\end{definition}

We call $R$ the anti-unitary antipode of $(A,\de)$.

\medskip

We know by \cite{JK} that $R_r \, \pi = \pi \, R_u$ which gives us, as usual, the next
result.

\begin{proposition}
We have for every $a \in A$ that $\pi_r(R(a)) = R_r(\pi_r(a))$.
\end{proposition}

\bigskip\medskip

In the rest of this section, we will prove the most basic relations between the objects
introduced in this section. In most cases, the proofs consist of pulling down the
corresponding relations on the \cst-algebra level. We will make use of the results on the
reduced \cst-algebra level because they were proved before the results on the universal
level (in fact, most of the latter results depend on the results on the reduced level).

\bigskip

First we describe the polar decomposition

\begin{result}
We have for every $z \in \C$ that $R \, \tau_z = \tau_z \, R$.
\end{result}

\begin{proposition} \label{art1.prop4}
We have that $S = R \, \tau_{-\frac{i}{2}}$.
\end{proposition}

These are immediate consequences of  corollary 5.4 of \cite{Kus} and
theorem 5.6 of \cite{Kus}.

\medskip

\begin{corollary} \label{art1.cor3}
We have the following commutation relations :
\begin{enumerate}
\item We have for every $z \in \C$ that $\tau_z \, S = S \, \tau_z$.
\item We have that $R \, S = S \, R$
\end{enumerate}
\end{corollary}

\bigskip

By corollary 7.3 of \cite{Kus}, we know that $(\si_r)_s \, (\tau_r)_t =
(\tau_r)_t \, (\si_r)_s$ for $s,t \in \R$. This will be used to prove the next result.

\begin{result}
We have for every $y,z \in \C$ that $\tau_y \, \si_z = \si_z \, \tau_y$.
\end{result}
\begin{demo}
The remarks before this result imply easily that $\si_t \, \tau_s = \tau_s \, \si_t$ for
$s,t \in \R$.

\vspace{1mm}

Take $a,b \in A$.

\vspace{1mm}

Fix $t \in \R$ for the moment. Then we have two analytic functions
$$\C \rightarrow \C : u \mapsto \vfi(\si_{-t}(b) \, \tau_u(a)) \hspace{1.5cm} \text{and}
\hspace{1.5cm} \C \rightarrow \C : u \mapsto \nu^u \, \vfi\bigl(\tau_{-u}(b) \,
\si_t(a)\bigr)$$
These functions are equal on the real line : we have for every $s \in \R$ that
$$\vfi(\si_{-t}(b) \, \tau_s(a))
= \vfi\bigl(b \, \si_t(\tau_s(a))\bigr) = \vfi\bigl(b \, \tau_s(\si_t(a))\bigr)
= \nu^s \, \vfi(\tau_{-s}(b) \, \si_t(a))$$
So they must be equal on $\C$. In particular $\vfi(\si_{-t}(b) \, \tau_y(a))
= \nu^y \, \vfi(\tau_{-y}(b) \, \si_t(a))$.  \ \ \ (*)

\medskip \vspace{1mm}

Now we have again two analytic functions
$$\C \rightarrow \C : u \mapsto \vfi(\si_{-u}(b) \, \tau_y(a)) \hspace{1.5cm}
\text{and} \hspace{1.5cm} \C \rightarrow \C : u \mapsto \nu^y \, \vfi(\tau_{-y}(b)
\, \si_u(a))$$
which by (*) are equal on the real line. So they must be equal on the whole complex plane.
In particular, $\vfi(\si_{-z}(b) \, \tau_y(a)) = \nu^y \, \vfi(\tau_{-y}(b) \, \si_z(a))$.
This gives us that
$$\vfi\bigl(b \, \tau_y(\si_z(a))\bigr) = \nu^y \, \vfi(\tau_{-y}(b) \, \si_z(a))
= \vfi(\si_{-z}(b) \, \tau_y(a)) = \vfi\bigl(b \, \si_y(\tau_z(a))\bigr)
$$
Hence the faithfulness of $\vfi$ implies that $\si_y \, \tau_z = \tau_z \, \si_y$.
\end{demo}

\bigskip

Now we prove some relations in connection with the comultiplication.

\begin{proposition} \label{art1.prop5}
Consider $z \in \C$. Then we have the following equalities :
\begin{itemize}
\item \ $(\tau_z \od \tau_z) \de = \de \tau_z$
\hspace{2.33cm} $\bullet$ \ $(\tau_z \od \si_z) \de = \de \si_z$
\item \ $(\si_z' \od \tau_{-z}) \de = \de \si_z'$
\hspace{2cm} $\bullet$ \ $(\si_z \od \si_{-z}') \de = \de \tau_z$
\end{itemize}
\end{proposition}
\begin{demo}
We only prove the first equality. The others are proven in the same way.

Choose $a \in A$. Take $b \in A$.

By proposition 5.7 of \cite{Kus}, we have for every $t \in \R$ that
\begin{eqnarray*}
& & (\pi_r \od \pi_r)(\de(\tau_t(a))(\tau_t(b) \ot 1))
= \de\bigl(\pi_r(\tau_t(a))\bigr)\,(\pi_r(\tau_t(b)) \ot 1)
= \de\bigl((\tau_r)_t(\pi_r(a))\bigr)\, ((\tau_r)_t(\pi_r(b)) \ot 1) \\
& & \spat = ((\tau_r)_t \ot (\tau_r)_t)\bigl(\de_r(\pi_r(a))\bigr)\,
((\tau_r)_t(\pi_r(b)) \ot 1)  = ((\tau_r)_t \ot (\tau_r)_t)(\de_r(\pi_r(a))\,
(\pi_r(b) \ot 1)) \\
& & \spat = ((\tau_r)_t \ot (\tau_r)_t)\bigl((\pi_r \od \pi_r)(\de(a)(b \ot 1))\bigr)
= (\pi_r \od \pi_r)\bigl((\tau_t \od \tau_t)(\de(a)(b \ot 1))\bigr)
\end{eqnarray*}
which implies that
$$\de(\tau_t(a))(\tau_t(b) \ot 1) = (\tau_t \od \tau_t)(\de(a)(b \ot 1))
\hspace{1.5cm} \text{(*)}$$

Choose $p,q \in A$. Then we have two functions
$$\C \rightarrow \C : y \mapsto (\vfi \od \vfi)((1 \ot p)\de(q) \de(\tau_y(a))
(\tau_y(b) \ot 1))$$
and
$$\C \rightarrow \C : y \mapsto (\vfi \od \vfi)\bigl((1 \ot p)\de(q)
(\tau_t \od \tau_t)(\de(a)(b \ot 1))\bigr)$$ We see immediately that the second function
is analytic. Because $(\vfi \od \vfi)((1 \ot p)\de(q) \de(\tau_y(a))(\tau_y(b) \ot 1)) =
\vfi(p \tau_y(b)) \, \vfi(q \tau_y(a))$ for $y \in A$, also the first function is
analytic.

We know by (*) that both functions are equal on the real axis, so they must be equal on
the whole complex plane. In particular, $(\vfi \od \vfi)((1 \ot p)\de(q)
\de(\tau_z(a))(\tau_z(b) \ot 1)) = (\vfi \od \vfi)\bigl((1 \ot p)\de(q) (\tau_z \od
\tau_z)(\de(a)(b \ot 1))\bigr)$.

Hence the faithfulness of $\vfi$ implies that
$$\de(\tau_z(a))(\tau_z(b) \ot 1) = (\tau_z \od \tau_z)(\de(a)(b \ot 1)) =
(\tau_z \od \tau_z)(\de(a)) (\tau_z(b) \ot 1)$$

\vspace{1mm}

So we see that $\de(\tau_z(a)) = (\tau_z \od \tau_z)(\de(a))$.

\end{demo}

\medskip

The next equality is proven in the same way as equality (*) in the first part of the
previous proof.

\begin{result}
We have that $\de R = \flip(R \od R)\de$.
\end{result}

\medskip

By the above equalities and the unicity of the co-unit on $(A,\de)$, we get the following
results.

\begin{corollary}
Consider $z \in \C$. Then $\vep \tau_z = \vep$.
\end{corollary}

\begin{corollary}
We have that $\vep R = \vep$.
\end{corollary}

\bigskip\medskip

In the next part, we look at some formulas involving the modular function $\sde$.

\medskip

The first one says that every $\sde^z$ is a  one-dimensional corepresentations of
$(A,\de)$.

\begin{proposition} \label{art1.prop6}
We have for every $z \in \C$ that $\de(\sde^z) = \sde^z \ot \sde^z$.
\end{proposition}
\begin{demo}
Take $a,b \in A$.

Then we have by proposition 8.6 of \cite{Kus} for every $t \in \R$ that
\begin{eqnarray*}
& & \hspace{-4mm}(\pi_r \od \pi_r)(\de(\sde^{it})\, \de(a)(b \ot 1)) =
(\pi_r \od \pi_r)(\de(\sde^{it}\,a)(b \ot 1))
= \de_r(\pi_r(\de(\sde^{it}\,a)) (\pi_r(b) \ot 1) \\
& & \hspace{-4mm} \spat = \de_r(\sde_r^{it} \, \pi_r(a))(\pi_r(b) \ot 1)
= (\sde_r^{it} \ot \sde_r^{it}) \, (\pi_r \od \pi_r)(\de(a) (b \ot 1))
= (\pi_r \od \pi_r)((\sde^{it} \ot \sde^{it}) \, \de(a)(b \ot 1))
\end{eqnarray*}
which implies that
$$\de(\sde^{it})\, \de(a)(b \ot 1) = (\sde^{it} \ot \sde^{it}) \, \de(a)(b \ot 1)
\hspace{1.5cm} (*)$$

We have now two functions
$$\C \rightarrow \C : u \mapsto (\vfi \od \vfi)(\de(\sde^{iu})\, \de(a)(b \ot 1))$$
and
$$\C \rightarrow \C : u \mapsto (\vfi \od \vfi)((\sde^{iu} \ot \sde^{iu})
\, \de(a)(b \ot 1))$$
The second function is clearly analytic. Because $(\vfi \od \vfi)(\de(\sde^{iu})\,
\de(a)(b \ot 1)) = \vfi(\sde^{iu} \, a) \, \vfi(b)$ for $u \in A$, also the second is
analytic. Furthermore, (*) implies that both functions are equal on the real line. So they
must be equal on the whole complex plane. In particular , $(\vfi \od \vfi)(\de(\sde^z)\,
\de(a)(b \ot 1))
= (\vfi \od \vfi)((\sde^z \ot \sde^z) \, \de(a)(b \ot 1))$

Hence the faithfulness of $\vfi$ implies that $\de(\sde^z)\, \de(a)(b \ot 1) = (\sde^z \ot
\sde^z) \, \de(a)(b \ot 1)$

\vspace{1mm}

So we get that $\de(\sde^z) = \sde^z \ot \sde^z$.
\end{demo}

\begin{corollary} \label{art1.cor2}
We have for every $z \in \C$ that $\vep(\sde^z) = 1$ and $S(\sde^z) = \sde^{-z}$.
\end{corollary}

\bigskip

\begin{result} \label{art1.res2}
Consider $y,z \in \C$. Then $\si_z(\sde^y) = \nu^{-y z} \, \sde^y$ and
$\tau_z(\sde^z) = \sde^{y z}$.
\end{result}
\begin{demo}
Choose $a \in A$.

Take $t \in \R$. By proposition 8.17 of \cite{Kus}, we know that $(\si_r)_t(\sde_r^y) =
\nu^{-y t} \, \sde_r^y$. This implies that
\begin{eqnarray*}
& & \pi_r(\si_t(\sde^y \, a)) = (\si_r)_t(\pi_r(\sde^y a)) = (\si_r)_t(\sde_r^y
\,\pi_r(a)) \\
& & \spat = (\si_r)_y(\sde_r^y) \, (\si_r)_t(\pi_r(a)) =  \nu^{-y t} \, \sde_r^y
\, \pi_r(\si_t(a)) = \nu^{-y t} \, \pi_r(\sde^y \, \si_t(a))
\end{eqnarray*}
which implies that $\si_t(\sde^y \, a) = \nu^{-y t} \, \sde^y \, \si_t(a)$. \ \ \ (*)

\medskip

Take $\om \in \ah$. Then we have two analytic functions
$$\C \rightarrow \C : c \mapsto \om(\si_c(\sde^y \, a))
\hspace{1cm} \text{and} \hspace{1cm}
\C \rightarrow \C : c \mapsto \nu^{-y c} \, \om(\sde^y \, \si_c(a))$$
which by (*) are equal on the real line. So they must be equal on the whole complex plane.
In particular, $\om(\si_z(\sde^y \, a)) = \nu^{-y z} \, \om(\sde^y \, \si_z(a))$.

Hence we get that $\si_z(\sde^y \, a) = \nu^{-y z} \, \sde^y \, \si_z(a)$, so
$\si_z(\sde^y) \, \si_z(a) = \nu^{-y z} \, \sde^y \, \si_z(a)$.

\medskip

So we see that $\si_z(\sde^y) = \nu^{-y z} \, \sde^y$.
The equality concerning $\tau$ is proven in the same way.
\end{demo}

\medskip

Combining the result concerning $\tau$ with corollary \ref{art1.cor2} and proposition
\ref{art1.prop4}, we get the following one.

\begin{result}
Consider $z \in \C$. Then $R(\sde^z) = \sde^{-z}$.
\end{result}

\bigskip

Definition 9.2 of \cite{Kus} implies for every $z \in \C$ that $(\si_r')_z = R_r \,
(\si_r)_{-z} \, R_r$. So we get the next result.

\begin{result}
Consider $z \in \C$. Then $\si_z' = R \, \si_{-z} \, R = S \, \si_{-z} \, S^{-1}
= S^{-1} \, \si_{-z} \, S$.
\end{result}

\medskip

Lemma \ref{art1.lem2} implies now the following formula for $\si'$.

\begin{result} \label{art1.res1}
We have for every $z \in \C$ and $a \in A$ that
$\si_z'(a) = \sde^{iz} \, \si_z(a) \, \sde^{-iz}$.
\end{result}

\bigskip

We end this section with some remarks concerning the right Haar functional on $(A,\de)$.
Recall that we have a right Haar functional $\vfi S$ on $(A,\de)$ but we do not know (yet)
whether $\vfi S$ is positive. By the formula $\flip(R \od R)\de = \de R$, we have however
the following proposition.

\begin{theorem}
The functional $\vfi R$ is a positive Haar functional on $(A,\de)$.
\end{theorem}

Because $S = \tau_{-\frac{i}{2}} \, R$, we have moreover that $\vfi S = \nu^{-\frac{i}{2}}
\, \vfi R$.

Because $\vfi S = \sde \vfi$, we have also that $\vfi(R(a)) = \vfi(\sde^\frac{1}{2} \, a
\, \sde^\frac{1}{2})$ for $a \in A$.

\bigskip

\section{The analytic structure of the dual}

In this section, we will connect the alytic objects associated to the dual of an algebraic
quantum groups to the analytic objects of this algebraic quantum groups.

\medskip

So consider an algebraic quantum group $(A,\de)$ with a positive left Haar functional
$\vfi$ on it. As in the previous section, we will use the notations $\si$ for the modular
group of the left Haar functional, $\si'$ for the modular group of the right Haar
functional, $\tau$  for the scaling group and $R$ for the anti-unitary antipode.

\medskip

The corresponding objects on the dual quantum group $(\ah,\deh)$ will get a hat on them,
e.g. $\hat{\si}$ will denote the modular group of the left Haar functional on
$(\ah,\deh)$.

\bigskip

First we start with the modular group of the dual quantum group $(\ah,\deh)$.
(Similar results are also considered in \cite{MasNak}).

\medskip

We will introduce first a temporary notation. Consider $\om \in M(\ah)$. Then we  define
$\om_z \in A'$ such that $\om_z(a) = \om(\tau_z(a) \, \sde^{-iz})$ for $a \in A$.

\medskip

\begin{lemma}
Consider $z \in \C$ and $a \in A$. Then
$$(\om_z \od \io)\de(a) = \tau_{-z}\bigl(\,[(\om \od \io)\de(\tau_z(a) \,
\sde^{-iz})] \,\, \sde^{iz}\,\bigr) \hspace{0.5cm} \text{and} \hspace{0.5cm}
(\io \od \om_z)\de(a) = \tau_{-z}\bigl(\,[(\io \od \om)\de(\tau_z(a)\,\sde^{-iz})]
 \,\, \sde^{iz}\,\bigr)$$
\end{lemma}
\begin{demo}
Take $b \in A$. Then we have that
$$[(\om_z \od \io)\de(a)] \, b = (\om_z \od \io)(\de(a)(1 \ot b))
= (\om \od \io)\bigl( (\tau_z \od \io)(\de(a) (1 \ot b))\, (\sde^{-iz} \ot 1) \bigr)
\hspace{1.5cm} \text{(*)}$$
Using proposition \ref{art1.prop5}, proposition \ref{art1.prop6} and result
\ref{art1.res2}, we see that
\begin{eqnarray*}
& & (\tau_z \od \io)(\de(a) (1 \ot b))\ (\sde^{-iz} \ot 1) =
(\io \od \tau_{-z})\bigl(\de(\tau_z(a)) \, (\sde^{-iz} \ot  \tau_z(b)) \bigr)  \\
& & \spat = (\io \od \tau_{-z}) \bigl(\de(\tau_z(a)\,\sde^{-iz}) \, (1 \ot \sde^{iz} \,
\tau_z(b))\bigr)
= (\io \od \tau_{-z}) \bigl(\de(\tau_z(a) \, \sde^{-iz}) \, (1 \ot \tau_z(\sde^{iz}
b))\bigr)
\end{eqnarray*}
So by applying $\om \od \io$ to this equation and using (*), we get that
\begin{eqnarray*}
& & [(\om_z \od \io)\de(a)] \, b
= \tau_{-z}\bigl((\om \od \io)(\de(\tau_z(a) \, \sde^{-iz}) \, (1 \ot \tau_z(\sde^{iz}
\, b)))\bigr)
=  \tau_{-z}\bigl([(\om \od \io)\de(\tau_z(a) \, \sde^{-iz})] \, \tau_z(\sde^{iz}
\, b)\bigl) \\
& & \spat =  \tau_{-z}\bigl([(\om \od \io)\de(\tau_z(a) \, \sde^{-iz})] \, \sde^{iz}
\, \tau_z(b)\bigl)
= \tau_{-z}\bigl([(\om \od \io)\de(\tau_z(a) \, \sde^{-iz})] \, \sde^{iz} \bigl) \  b
\end{eqnarray*}
where we used result \ref{art1.res2} in the second last equality.

Consequently, $(\om_z \od \io)\de(a) =  \tau_{-z}\bigl([(\om \od \io)\de(\tau_z(a)
\sde^{-iz})] \, \sde^{iz} \bigl)$. The other equality is proven in the same way.
\end{demo}

\medskip

Using theorem 2.12 of \cite{JK3}, it is now easy to verify the next lemma.

\begin{lemma} \label{art4.lem1}
Consider $z \in \C$. Then we have the following properties.
\begin{enumerate}
\item We have for every $\om \in M(\ah)$ that $\om_z \in M(\ah)$.
\item Consider $\om,\th \in M(\ah)$. Then $(\om \th)_z = \om_z \, \th_z$.
\end{enumerate}
\end{lemma}

\medskip

\begin{lemma} \label{art4.lem2}
Consider $z \in \C$ and $\om \in \ah$. Then $\om_z$ belongs to $\ah$ and
$\hat{\vfi}(\om_z) = \hat{\vfi}(\om)$
\end{lemma}
\begin{demo}
There exist $a \in A$ such that $\om = \psi a$.

Result \ref{art1.res2} and corollary \ref{art1.cor3} imply that $\psi \, \tau_z = \nu^z \,
\psi$. Results \ref{art1.res1} and \ref{art1.res2} imply that $\si_i'(\sde^{-iz}) =
\nu^{-z} \, \sde^{-iz}$.

So we have for every $x \in A$ that
\begin{eqnarray*}
\om_z(x) & = & (\psi a)(\tau_z(x)  \, \sde^{-iz}) = \psi(a \, \tau_z(x) \, \sde^{-iz})
= \psi\bigl(\tau_z(\tau_{-z}(a)\, x \, \sde^{-iz})) \\
& = & \nu^z \, \psi(\tau_{-z}(a) \, x \, \sde^{-iz}) = \psi(\sde^{-iz} \, \tau_{-z}(a) \, x )
\end{eqnarray*}
So we see $\om_z = \psi \sde^{-iz} \tau_{-z}(a)$. This implies that $\om_z \in \ah$ and
that $$\hat{\vfi}(\om_z) = \vep(\sde^{-iz} \tau_{-z}(a)) = \vep(a)
= \hat{\vfi}(\om)$$
\end{demo}

\medskip

\begin{lemma} \label{art4.lem3}
Consider $\om \in \ah$ and $a \in A$. Then the fuction $\C \rightarrow \C : z \mapsto
\om_z(a)$ is analytic.
\end{lemma}
\begin{demo}There exist $b,c,d \in A$ such that $\om = b c \vfi d^*$.

We know that $\pi_r(a)$ is analytic with respect to $\tau_r$ and that $\pi_r(b)$ is
analytic with respect to $\sde_r$. This implies that the function $\C \rightarrow A_r : z
\mapsto (\tau_r)_z(\pi_r(a)) \, \sde_r^{-iz} \, \pi_r(b)$ is analytic. So the function $\C
\rightarrow A_r : z \mapsto \vfi_r(\pi_r(d)^* \, (\tau_r)_z(\pi_r(a)) \, \sde_r^{-iz} \,
\pi_r(b) \, \pi_r(c))$ is analytic.

But we have for every $z \in \C$ that
\begin{eqnarray*}
& & \vfi_r(\pi_r(d)^* \, (\tau_r)_z(\pi_r(a)) \, \sde_r^{-iz} \, \pi_r(b) \, \pi_r(c))
= \vfi_r(\pi_r(d)^* \pi_r(\tau_z(a)) \pi_r(\sde^{-iz} \, b) \pi_r(c)) \\
& & \spat
= \vfi_r(\pi_r(d^* \tau_z(a) \sde^{-iz} b c)) = \vfi(d^* \tau_z(a) \sde^{-iz} b c)
= \om(\al_z(a))
\end{eqnarray*}

Hence the function  $\C \rightarrow \C : z \mapsto \om(\al_z(a))$ is analytic.
\end{demo}

\bigskip

\begin{proposition}
We have for every $\om \in \ah$, $a \in A$ and $z \in \C$ that
$\hat{\si}_z(\om)(a) = \om(\tau_z(a) \, \sde^{-iz})$.
\end{proposition}
\begin{demo}
Define the function $\be$ from $\C$ into $L(\ah)$ such that $\be_z(\om) = \om_z$ for $z
\in \C$ and $\om \in \ah$. Then we have the following properties.
\begin{enumerate}
\item It is clear that $\be_0 = \io$
\item Choose $s,t \in \R$.

Take $\om \in \ah$. Then we have for every $a \in A$ that
\begin{eqnarray*}
& & [\be_s(\be_t(\om))](a) =  [\be_t(\om)](\tau_s(a) \, \sde^{-it})
= \om(\tau_t(\tau_s(a) \, \sde^{-it}) \, \sde^{-is}) \\
& & \spat =  \om(\tau_t(\tau_s(a)) \ \sde^{-it} \, \sde^{-is})
= \om(\tau_{s+t}(a) \,  \sde^{-i(s+t)})) = [\be_{s+t}(\om)](a) \ ,
\end{eqnarray*}
implying that $\be_s(\be_t(\om)) = \be_{s+t}(\om)$.

\vspace{1mm}

So we have proven that $\be_s \, \be_t = \be_{s+t}$.

\item Lemma \ref{art4.lem1} implies that $\be_t$ is multiplicative for every $t \in \R$.

\item Choose $t \in \R$.

Take $\om \in \ah$. Then we have for every $a \in A$ that
\begin{eqnarray*}
& & [\be_t(\om)]^*(a)  =  \overline{[\be_t(\om)](S(a)^*)}
= \overline{\om(\tau_t(S(a)^*) \, \sde^{-it})}
=  \overline{\om((\sde^{it} \, \tau_t(S(a)))^*)} \\
&  & \spat = \overline{\om((\sde^{it} \, S(\tau_t(a)))^*)}
= \overline{\om(S(\tau_t(a) \,\sde^{-it})^*)}
= \om^*(\tau_t(a) \,\sde^{it}) = [\be_t(\om^*)](a) \ ,
\end{eqnarray*}
implying that $[\be_t(\om)]^* = \be_t(\om^*)$.

\vspace{1mm}

So we have proven that $\be_t$ preserves the $^*$-operation.

\item Lemma \ref{art4.lem2} implies that $\hat{\vfi}$ is invariant under $\be_t$
for every $t \in \R$.

\item Choose $b \in \ahh$ and $\om \in \ah$. Then there exist $a \in A$ such that
$b(\th) = \th(a)$ for $\th \in \ah$. So we have for every $z \in \C$ that
$b(\be_z(\om)) = \be_z(\om)(a) = \om_z(a))$ which implies that the function
$\C \rightarrow \C : z \mapsto  b(\be_z(\om))$ is analytic by lemma \ref{art4.lem3}
\end{enumerate}
So we can conclude from this al that $\be_z$ is an analytic one-parameter group on $\ah$.

\medskip

Choose $a \in A$. Lemma 2.8 of \cite{JK3} implies that $\hat{\si}_i(\psi a) = \psi \sde
S^2(a)$. We know by the proof of lemma \ref{art4.lem2} that $(\psi a)_i = \psi \sde
S^2(a)$. So we see that $\hat{\si}_i(\psi a) =  \be_i(\psi a)$.

Hence $\hat{\si}_i = \be_i$ which by proposition \ref{art2.prop1} implies that $\hat{\si}
= \be$.
\end{demo}

\medskip

The next result follows now immediately from lemma \ref{art4.lem1}.

\begin{corollary}
We have for every $\om \in M(\ah)$, $a \in A$ and $z \in \C$ that
$\hat{\si}_z(\om)(a) = \om(\tau_z(a) \, \sde^{-iz})$.
\end{corollary}

\bigskip

The proof of the next result is completely similar (an easier) to the proof of the
previous proposition. It is a consequence of the fact that $\hat{S}^2(\om) = \om S^2$ for
$\om \in \ah$.

\begin{proposition}
We have for every $\om \in M(\ah)$, $a \in A$ and $z \in \C$ that
$\hat{\tau}_z(\om) = \om \tau_z$.
\end{proposition}

\medskip

Combining this with the fact that $\hat{R} = \hat{S} \, \hat{\tau}_{\frac{i}{2}}$, we get
easily the following on.

\begin{corollary}
We have for every $\om \in M(\ah)$ that $\hat{R}(\om) = \om R$.
\end{corollary}

\medskip

Remembering that $\hat{\si}_z' = \hat{R} \, \hat{\si}_{-z} \, \hat{R}$, it is now easy to
check the next equality.

\begin{corollary}
We have for every $a \in A$, $\om \in M(\ah)$ and $z \in \C$ that
$(\hat{\si}_z')(\om)(a) = \om\bigl(\sde^{-iz} \, \tau_{-z}(a))\bigr)$
\end{corollary}

\bigskip\medskip

In the last part of this section, we look a little bit further into the modular function
$\hat{\sde}$ of the dual and connect it with the modular group $\si$ and the antipode
$\vep$.

\medskip

\begin{lemma}
We have that $\vep \si_z = \vep \si_z'$ for $z \in \C$.
\end{lemma}

This follows from the fact that $\si_z'(a) = \sde^{iz} \, \si_z(a) \, \sde^{-iz}$ for $a
\in A$ and the fact that $\vep(\sde^{iz}) = 1$.

\begin{result}
Consider $z \in \C$. Then $\vep \si_z$ belongs to $M(\ah)$.
\end{result}
\begin{demo}
Choose $a \in A$. By using proposition \ref{art1.prop5}, we see that
$$(\io \od \vep \si_z)\de(a) = (\tau_{-z} \od \vep)\de(\si_z(a)) = \tau_{-z}(\si_z(a))$$
So we get that $(\io \od \vep \si_z)\de(a) \in A$.

Using the fact that $\vep \si_z = \vep \si_z'$ and proposition \ref{art1.prop5} once
again, we get in a similar way that $(\vep \si_z \od \io)\de(a)$ belongs to $A$. Hence
theorem 2.12 of \cite{JK3} implies that $\vep \si_z \in M(\ah)$.
\end{demo}

\medskip

We introduce again a temporary notation. Define the mapping $u$ from $\C$ into $M(\ah)$
such that $u_z = \vep \si_{-z}$ for every $z \in \C$.

\vspace{1mm}

Recall from the proof of the previous lemma
that $(\io \od u_z)\de(a) = \tau_z(\si_{-z}(a))$ for $z \in \C$ and $a \in A$.

\medskip

\begin{lemma}
The mapping $u$ is an analytic unitary representation on $\ah$.
\end{lemma}
\begin{demo}
\begin{enumerate}
\item It is clear that $u_0 = 1$.
\item Choose $s,t \in A$. Take $a \in A$. By the proof of the previous result,
we know that $(\io \od u_t)\de(a) = \tau_t(\si_{-t}(a))$. This implies by
theorem 2.12 of \cite{JK3} that
\begin{eqnarray*}
(u_s \, u_t)(a) & = & u_s((\io \od u_t)\de(a)) =
\vep\bigl(\si_{-s}(\tau_t(\si_{-t}(a)))\bigr) \\ & = &
\vep\bigl(\tau_t(\si_{-s}(\si_{-t}(a)))\bigr) = \vep(\si_{-(s+t)}(a)) = u_{s+t}(a)
\end{eqnarray*}
So we see that $u_{s+t} = u_s \, u_t$.
\item Choose $t \in \R$. Take $a \in A$.
Then $$u_t^*(a) = \overline{u_t(S(a)^*)} = \overline{\vep(\si_{-t}(S(a)^*))}
= \vep(\si_{-t}(S(a))) = \vep(S(\si_t'(a))) = \vep(\si_t'(a)) = u_{-t}(a)$$
So we get that $u_t^* = u_{-t}$.
\item Choose $a \in \ahh$. Then there exists $\om \in \ah$ such that $a = \om \, a$.
We know that their exists $b \in A$ such that $a(\th) = \th(b)$ for every $\th \in \ah$.
Define $c = (\io \ot \om)\de(b) \in A$. Recall that $\pi_u(c)$ is analytic with respect to
$\si_u$.

We have moreover for every $z \in \C$ that
$$a(u_z) = a(u_z \om) = (u_z \om)(b) = u_z((\io \ot \om)\de(b)) = u_z(c)
= \vep_u\bigl((\si_u)_z(\pi_u(c))\bigr)$$ which implies that
the function $\C \rightarrow \C : z \mapsto a(u_z)$ is analytic.
\end{enumerate}
\end{demo}

\medskip

\begin{lemma}
We have that $u_{-i} = \hat{\sde}$.
\end{lemma}
\begin{demo}
Take $\th \in \ah$. Then there exist $a \in A$ such that $\th = \psi a$.

So we have for every $x \in A$ that
\begin{eqnarray*}
(\th S \rho^{-1})(x) & = & \psi\bigl(a\,S(\rho^{-1}(x))\bigr) = \psi\bigl(a
\,\rho'(S(x))\bigr)
= \psi(S(x) \, a) \\
& = & \vfi(S(x)\,a\,\sde) = \vfi(S( \sde^{-1} S^{-1}(a) \, x)) = \psi(\sde^{-1}
S^{-1}(a) \, x)
\end{eqnarray*}
So we see that $\th S \rho^{-1} = \psi \sde^{-1} S^{-1}(a)$ which implies that
$$\hat{\vfi}(\th S \rho^{-1}) = \vep(\sde^{-1} S^{-1}(a)) = \vep(a) =
\hat{\vfi}(\th) \hspace{1.5cm} \text{(*)}$$

\medskip

Choose $\om \in A$. Then we have for every $x \in A$ that
$$(\om u_{-i})(x) = \om((\io \od u_{-i})\de(x)) = \om\bigl(\tau_{-i}(\si_i(x))\bigr)
= \om\bigl(S^2(\rho^{-1}(x))\bigr)$$
so we get that $\om u_{-i} = \om S^2 \rho^{-1} = (\om S) S \rho^{-1}$. Therefore equation
(*) gives us that
$$\hat{\vfi}(\om u_{-i}) = \hat{\vfi}((\om S) S \rho^{-1}) = \hat{\vfi}(\om S)
= \hat{\vfi}(\hat{S}(\om))$$

This implies that $u_{-i} = \hat{\sde}$.
\end{demo}

\medskip

Proposition \ref{art1.prop7} implies  now immediately the next proposition.

\begin{proposition}
We have for every $z \in \C$ that $\hat{\sde}^{i z} = \vep \si_{-z}$.
\end{proposition}

\bigskip\medskip

In the last part of this section, we will prove a connection between $\hat{\sde}$ and the
modular functionals introduced by Woronowicz in \cite{Wor1}.

\medskip

\begin{lemma}
Consider $a \in A$ and $z \in \C$.
Then \begin{enumerate}
\item $(\vep \si_z \od \io \od \vep \si_{-z})\de^{(2)}(a) = \sde^{iz} \, \tau_{2z}(a)
\, \sde^{-iz}$
\item $(\vep \si_z \od \io \od \vep \si_z)\de^{(2)}(a) = \sde^{iz} \, \si_{2z}(a)
\, \sde^{-iz}$
\end{enumerate}
\end{lemma}
\begin{demo}
In this proof, we will make use of proposition \ref{art1.prop5} and result
\ref{art1.res1}.
\begin{enumerate}
\item We have that
\begin{eqnarray*}
& & (\vep \si_z \od \io \od \vep \si_{-z})\de^{(2)}(a) = (\vep \si_z \od \io \od
\vep \si_{-z})((\io \od \de)\de(a)) = (\vep \si_z \od \tau_z \si_{-z}) \de(a) \\
& & \spat = \sde^{iz} \, \, (\vep \si_z \od  \tau_z \si_{-z}') \de(a) \, \, \sde^{-iz}
= \sde^{iz} \, \tau_{2z}(a)  \, \sde^{-iz}
\end{eqnarray*}
\item And similary,
\begin{eqnarray*}
& & (\vep \si_z \od \io \od \vep \si_z)\de^{(2)}(a) = (\vep \si_z \od \io \od
\vep \si_z)((\io \od \de)\de(a)) = (\vep \si_z \od \tau_{-z} \si_z) \de(a) \\
& & \spat = (\vep \si_z' \od \si_z \tau_{-z}) \de(a) = \si_z(\si_z'(a))
= \sde^{iz} \, \si_{2 z}(a) \, \sde^{-iz}
\end{eqnarray*}
\end{enumerate}
\end{demo}

\medskip

In theorem 5.6 of \cite{Wor1}, Woronowicz introduced for a compact quantum group the
linear functionals  $f_z  \in A'$ for every complex number $z$. The above lemma guarantees
that they are connected with the modular function $\hat{\sde}$.

\medskip

\begin{result}
Suppose that $(A,\de)$ is a compact algebraic quantum group. Then we have for every $z \in
\C$ that $\sdeh^z = f_{-2z}$.
\end{result}
\begin{demo} Recall that $\sde = 1$ in this case.

Because $A$ has a unit, theorem 2.12 of \cite{JK3} implies immediately that $M(\ah) = A'$.
So we get for every $z \in \C$ that $f_z$ belongs to $M(\ah)$. Define the function $v$
from $\C$ into $M(\ah)$ such that $v_z = f_{- 2 i z}$ for $z \in \C$. Then theorem 5.6 of
\cite{Wor1} implies immediately that $v$ satisfies the following properties :
\begin{enumerate}
\item $v_0 = 1$
\item We have for every $s,t \in \C$ that $v_{s+t} = v_s \, v_t$.
\item We have for every $t \in \R$ that $v_t^* = v_{-t}$.
\item Consider $\om \in \ahh$. Then the function $\C \rightarrow \C : z \rightarrow
\om(v_z)$ is analytic (remember that $\ahh = A$).
\end{enumerate}
So $v$ is an analytic unitary representation on $\ah$.

\vspace{1mm}

Theorem 5.6 of \cite{Wor1} implies also that $f_1$ is invertible in $M(\ah)$ and that
$(f_1)^{-1} =f_{-1}$.

Choose $\om \in \ah$.

Take $a \in A$. Then item 5 of theorem 5.6 of \cite{Wor1} and the previous lemma imply
that
$$(f_1 \od \io \od f_{-1})\de^{(2)}(a) = S^{2}(a) = \tau_{-i}(a) = (\sdeh^{-\frac{1}{2}}
\od \io \od \sdeh^{\frac{1}{2}})\de^{(2)}(a)$$
So by applying $\om$ to this equality, we get that
$(f_1 \, \om \, f_{-1})(a) = (\sdeh^{-\frac{1}{2}} \, \om \, \sdeh^{\frac{1}{2}})(a)$.

Hence $f_1 \, \om \, f_{-1} = \sdeh^{-\frac{1}{2}} \, \om \, \sdeh^{\frac{1}{2}}$, which
implies that $\sdeh^{\frac{1}{2}} \, f_1 \, \om = \om \, \sdeh^{\frac{1}{2}} \, f_1$. \ \
\ (a)

\vspace{1mm}

Item 6 of theorem 5.6 of \cite{Wor1} implies that $f_1 * a * f_1 = \rho(a)$ for $a \in A$.
Combining this with the previous lemma, we get in a similar way that $f_1 \, \om \, f_1 =
\sdeh^{-\frac{1}{2}} \, \om \, \sdeh^{-\frac{1}{2}}$ which implies that
$\sdeh^{\frac{1}{2}} \, f_1 \, \om = \om \, \sdeh^{-\frac{1}{2}} \, f_{-1}$. \ \ \ (b)

\vspace{1mm}

Combining (a) and (b), we get that $\om \, \sdeh^{\frac{1}{2}} \, f_1 = \om
\,\sdeh^{-\frac{1}{2}} \, f_{-1}$.

\medskip

So we see that $\sdeh^{\frac{1}{2}} \, f_1 = \sdeh^{-\frac{1}{2}} \, f_{-1}$ which implies
that $\sdeh = f_{-2}$. Consequently, $v_{-i} = \sdeh$.

Therefore proposition \ref{art1.prop7} implies for every $z \in \C$ that $\sdeh^z =
v_{-iz} = f_{-2 z}$.
\end{demo}

\end{document}